\documentclass[a4paper,reqno]{amsart}
\usepackage{amssymb,upref,mathrsfs,url}
\let\mathcal\mathscr
\usepackage[mathscr]{eucal}

\usepackage[matrix,arrow,tips,ps,dvips]{xy}

\SelectTips{cm}{}

\allowdisplaybreaks

\newcommand*{\pd}[2]{\mathchoice{\frac{\partial#1}{\partial#2}}
 {\partial#1/\partial#2}{\partial#1/\partial#2}
 {\partial#1/\partial#2}}
\newcommand*{\fd}[2]{\mathchoice{\frac{\delta#1}{\delta#2}}
 {\delta #1/\delta#2}{\delta#1/\delta#2}{\delta#1/\delta#2}}

\newcommand{\eval}[2][\right]{\relax
 \ifx#1\right\relax \left.\fi#2#1\rvert}

\newcommand{\envert}[2][\right]{\relax
 \ifx#1\right\relax \left\lvert\else#1\lvert\fi#2#1\rvert}
\let\abs=\envert

\newcommand{\enVert}[2][\right]{\relax
 \ifx#1\right\relax \left\lVert\else#1\lVert\fi#2#1\rVert}
\let\matr=\enVert

\let\kappa\varkappa
\let\phi\varphi

\DeclareMathOperator{\Hom}{Hom}
\DeclareMathOperator{\im}{im}

\DeclareMathOperator{\CDiff}{\mathcal{C}Dif{}f}

\newcommand*{\CDiffskew}[1]{\CDiff_{(#1)}^{\,\text{\textup{skew}}}}
\newcommand*{\CDiffskewad}[1]{\CDiff_{(#1)}^{\,\text{\textup{sk-ad}}}}

\newcommand{\cprime}{\/{\mathsurround=0pt$'$}}

\providecommand{\href}[2]{#2}
\newcommand{\urlprefix}{URL }
\newcommand*{\eprint}[2][]{%
\href{http://arXiv.org/abs/#2}{\begingroup \Url{arXiv:#2}}%
}

\DeclareFontFamily{OML}{cyi}{}
\DeclareFontShape{OML}{cyi}{m}{n}{
  <5> <6> <7> <8> <9> gen * wncyi
  <10> <10.95> <12> <14.4> <17.28> <20.74> <24.88> wncyi10
 }{}
\DeclareSymbolFont{rusletters}{OML}{cyi}{m}{n}
\DeclareSymbolFontAlphabet{\rusmath}{rusletters}
\DeclareMathSymbol\re{\rusmath}{rusletters}{"03}

\newcommand{\Eu}{\mathscr{E}}

\newtheorem{theorem}{Theorem}
\newtheorem{proposition}{Proposition}

\theoremstyle{definition}
\newtheorem{definition}{Definition}
\newtheorem{example}{Example}

\theoremstyle{remark}
\newtheorem{remark}{Remark}

\begin{document}

\hfill DIPS-6/2002

\hfill math.DG/0304245

\vspace{1cm}

\title{Hamiltonian operators and $\ell^*$-coverings}

\author{P. Kersten}

\address{Paul Kersten \\
  University of Twente, Faculty of Electrical Engineering, Mathematics
  and Computer Science \\
  P.O.~Box 217 \\
  7500 AE Enschede \\
  The Netherlands}

\email{kersten@math.utwente.nl}

\author{I. Krasil{\cprime}shchik}

\address{Iosif Krasil{\cprime}shchik \\
  The Diffiety Institute and Independent University of Moscow \\
  B.~Vlasevsky~11 \\
  119002 Moscow \\
  Russia}

\email{josephk@diffiety.ac.ru}

\author{A.~Verbovetsky}

\address{Alexander Verbovetsky \\
  Independent University of Moscow \\
  B. Vlasevsky~11 \\
  119002 Moscow \\
  Russia}

\email{verbovet@mccme.ru}


\begin{abstract}
  An efficient method to construct Hamiltonian structures for
  nonlinear evolution equations is described. It is based on the
  notions of variational Schouten bracket and $\ell^*$-covering.  The
  latter serves the role of the cotangent bundle in the category of
  nonlinear evolution PDEs.  We first consider two illustrative
  examples (the KdV equation and the Boussinesq system) and
  reconstruct for them the known Hamiltonian structures by our
  methods.  For the coupled KdV-mKdV system, a new Hamiltonian
  structure is found and its uniqueness (in the class of polynomial
  $(x,t)$-independent structures) is proved. We also construct a
  nonlocal Hamiltonian structure for this system and prove its
  compatibility with the local one.
\end{abstract}
\maketitle

\section*{Introduction}
We describe a method of constructing Hamiltonian structures for
nonlinear evolution equations (or systems of such equations). The
method is based on two concepts: the \emph{variational Schouten
bracket} and the \emph{$\ell^*$-covering} over a nonlinear PDE.

In Section~\ref{sec:gener-jet-bundl}, we expose some general facts
concerning the geometry of super PDE. In
Section~\ref{sec:vari-scho-brack}, we construct the variational
Schouten bracket on a super version of Kupershmidt's cotangent bundle
to a bundle and, following~\cite{IgoninVerbovetskyVitolo:FLVDOp},
obtain an explicit formula for this bracket. In
Section~\ref{sec:hamilt-evol-equat}, simple computational formulas are
deduced to check the Hamiltonianity of a bivector and compatibility of
two Hamiltonian structures.  Using the Schouten bracket, we define
Hamiltonian evolution equations (including the cases when the
Hamiltonian operator~$A$ may depend explicitly on time while the
equation itself may not possess a Hamiltonian functional).  This
definition is equivalent to the operator equality
\begin{equation*}
  \ell_\mathcal{E}\circ A+A\circ\ell^*_\mathcal{E}=0,
\end{equation*}
where $\ell_\mathcal{E}$ is the linearization of the equation and
$\ell^*_\mathcal{E}$ is the adjoint operator. To solve this equation,
we introduce the notion of $\ell^*_\mathcal{E}$-covering (which is a
particular case of a more general construction introduced in
Section~\ref{sec:delta-coverings}) and show that to any operator $A$
satisfying the above equation there corresponds a function $s$ on the
$\ell^*_\mathcal{E}$-covering such that
$\tilde{\ell}_\mathcal{E}(s)=0$, where $\tilde{\ell}_\mathcal{E}$ is
the lifting of the linearization operator for $\mathcal{E}$ to the
$\ell^*_\mathcal{E}$-covering.  In other words, the operators we are
interested in are identified with shadows of nonlocal symmetries (of a
special type) in the $\ell^*_\mathcal{E}$-covering.

The reason to introduce the concept of $\ell^*$-covering is two-fold.
First, as it was just indicated, it allows to reduce construction of
Hamiltonian structures to computation of symmetries (with a subsequent
check of additional conditions) for which a number of efficient
software packages exists. Second, to our opinion, this point of view
gives a new and fruitful insight into the theory of Hamiltonian
structures for partial differential equations.

In Sections~\ref{sec:korteweg-de-vries}
and~\ref{sec:hamilt-struct-bouss}, these methods are applied to the
known examples of the KdV equation and the Boussinesq system.  In
Section~\ref{sec:appl}, we construct a Hamiltonian structure for the
coupled KdV-mKdV system~\cite{KerstenKrasilshchik:CInCKS}. We also
prove that this structure is unique in the class of
$(x,t)$-independent polynomial structures.  Nevertheless, extending
the initial setting with certain nonlocal variables, we find another
Hamiltonian operator that serves a Hamiltonian structure for `higher'
coupled KdV-mKdV equations. This structure is compatible with the
local one. It is to be noted that the theory of nonlocal Hamiltonian
structures is not sufficiently developed yet and needs additional
research.

In the Appendix, we briefly recall the construction of the recursion
operator for the coupled KdV-mKdV system (obtained earlier
in~\cite{KerstenKrasilshchik:CInCKS}) by which the above mentioned
Hamiltonian structures are related to each other.

\section{Generalities: Jet bundles and differential equations}
\label{sec:gener-jet-bundl}

Let us formulate the main definitions and results we will use.  For
more details we refer
to~\cite{KrasilshchikVinogradov:SCLDEqMP,KrasilshchikVerbovetsky:HMEqMP,%
IgoninVerbovetskyVitolo:FLVDOp}.

\subsection{Jet bundles}
\label{sec:jet-bundles}

Let $\pi\colon E\to M$ be a vector bundle over an $n$-dimensional base
manifold~$M$ and $\pi_\infty\colon J^\infty(\pi)\to M$ be the infinite
jet bundle of local sections of the bundle~$\pi$.

In coordinate language, if $x_1,\dots,x_n$,~$u^1,\dots,u^m$ are
coordinates on~$E$ such that $x_i$ are base coordinates and $u^j$ are
fiber ones, then $\pi_\infty\colon J^\infty(\pi)\to M$ is an
infinite-dimensional vector bundle with fiber coordinates $u^j_\tau$,
where $\tau=i_1\dots i_{\abs{\tau}}$ is a symmetric multi-index.

Now, we generalize the definition of the jet bundle to the case of
\emph{superbundles}.

\begin{definition}
  Let $E$ be a supermanifold of superdimension $(n+m_0)\vert m_1$, and
  $\pi\colon E\to M$ be a vector bundle over an $n$-dimensional even
  manifold~$M$.  If $\pi$ is split into the direct sum of two vector
  subbundles $\pi=\pi^0\oplus\pi^1$ such that the fibers of $\pi^0$
  are even and the fibers of $\pi^1$ are odd, then we say that $\pi$
  (along with the splitting) is a \emph{superbundle}.
\end{definition}

For a superbundle~$\pi$, we define the \emph{infinite jet
superbundle~$\pi_\infty\colon J^\infty(\pi)\to M$} by setting:
\begin{align*}
  (\pi_\infty)^0&=(\pi^0)_\infty, \\
  (\pi_\infty)^1&=(((\pi^1)^\Pi)_\infty)^\Pi,
\end{align*}
where the superscript ${}^\Pi$ denotes the reversion of parity.

Denote by $\mathcal{F}(\pi)$ the superalgebra of smooth functions
on~$J^\infty(\pi)$.

\begin{remark}
  By definition, we have
  \begin{equation*}
    \mathcal{F}(\pi)=\mathcal{F}(\pi^0)
    \otimes_{C^\infty(M)}
    \Lambda^*(\mathcal{F}_{\text{\textup{lin}}}((\pi^1)^\Pi)),
  \end{equation*}
  where $\mathcal{F}_{\text{\textup{lin}}}(\,\cdot\,)
  \subset\mathcal{F}(\,\cdot\,)$ is the subspace of functions linear
  along fibers.
\end{remark}

In what follows we shall use the term `bundle' to mean `superbundle'.

\subsection{The Cartan distribution}
\label{sec:cartan-distribution}

Consider a bundle~$\pi\colon E\to M$ and define the
\emph{$C^\infty(M)$-supermodule $\Gamma(\pi)$} of its `sections' as
follows.  If $\pi$ is even, then $\Gamma(\pi)$ is the module of
sections of~$\pi$.  If $\pi$ is a general superbundle, then we put
$\Gamma(\pi)=\Gamma(\pi)^0\oplus\Gamma(\pi)^1$, with
$\Gamma(\pi)^0=\Gamma(\pi^0)$ and
$\Gamma(\pi)^1=(\Gamma((\pi^1)^\Pi))^\Pi$.

\begin{remark}
  Thus, in line with our definition of jets of superbundles, we define
  elements of $\Gamma(\pi)$ to be pairs of sections of $\pi^0$ and
  $\pi^1$.
\end{remark}

Next, we note that every fiberwise linear function~$f$ on infinite
jets $J^\infty(\pi)$ can be naturally identified with a linear
differential operator~$\nabla_f\colon\Gamma(\pi)\to C^\infty(M)$ and
vice versa.  Indeed, for even bundle~$\pi$ the correspondence is given
by the relation
\begin{equation*}
  \nabla_f(s)(a)=f(j_\infty(s)(a)),
\end{equation*}
where $s\in\Gamma(\pi)$, $j_\infty(s)$ is the infinite jet of~$s$,
$a\in M$.  The general case reduces to the even one, since
$\mathcal{F}_{\text{\textup{lin}}}(\pi)
=\mathcal{F}_{\text{\textup{lin}}}(\pi^0\oplus(\pi^1)^\Pi)$.

\begin{remark}
  \label{sec:cartan-distribution-1}
  The maps
  \begin{align*}
    j_\infty&\colon\Gamma(\pi^0)\to\Gamma(\pi_\infty^0), \\
    j_\infty&\colon\Gamma((\pi^1)^\Pi)\to\Gamma((\pi_\infty^1)^\Pi)
  \end{align*}
  give rise to a map of supermodules
  $j_\infty\colon\Gamma(\pi)\to\Gamma(\pi_\infty)$.
\end{remark}

The infinite jet bundle $\pi_\infty\colon J^\infty(\pi)\to M$ admits a
natural flat connection such that the lift~$\hat X$ of a vector field
$X$ on~$M$ is uniquely defined by the condition
\begin{equation*}
  \nabla_{\hat X(f)}=X\circ\nabla_f,\quad f\in\mathcal{F}_{\text{\textup{lin}}}(\pi).
\end{equation*}
In coordinates, the lift of $\pd{}{x_i}$ is the $i$th \emph{total
derivative}
\begin{equation*}
  \frac{\widehat\partial}{\partial x_i}=D_i
  =\pd{}{x_i}+\sum_{j,\tau}u^j_{\tau i}\pd{}{u^j_\tau}.
\end{equation*}
Vector fields of the form~$\hat X$ generate an $n$-dimensional
distribution on~$J^\infty(\pi)$ called the \emph{Cartan distribution}
and denoted by~$\mathcal{C}(\pi)$.  Obviously, the Cartan distribution
is Frobenious in the sense that
$[\mathcal{C}(\pi),\mathcal{C}(\pi)]\subset\mathcal{C}(\pi)$.  In
coordinate language, the Cartan distribution is spanned by the total
derivatives.

\subsection{Horizontal calculus and evolutionary fields}
\label{sec:horiz-calc-evol}

Let $\xi\colon B\to M$ be a vector bundle and $\pi_\infty^*(\xi)\colon
B\times_M J^\infty(\pi)\to J^\infty(\pi)$ its pullback
along~$\pi_\infty$.  The $C^\infty(J^\infty(\pi))$-supermodule
$\Gamma(\pi_\infty^*(\xi))$ is defined as above:
$\Gamma(\pi_\infty^*(\xi))=\Gamma(\xi)
\otimes_{C^\infty(M)}C^\infty(J^\infty(\pi))$, if $\xi$ is even, and
$\Gamma(\pi_\infty^*(\xi))
=\Gamma(\pi_\infty^*(\xi))^0\oplus\Gamma(\pi_\infty^*(\xi))^1$, with
$\Gamma(\pi_\infty^*(\xi))^0=\Gamma(\pi_\infty^*(\xi)^0)$ and
$\Gamma(\pi_\infty^*(\xi))^1=(\Gamma((\pi_\infty^*(\xi)^1)^\Pi))^\Pi$
if $\xi$ is a general superbundle.
\begin{definition}
  A $C^\infty(J^\infty(\pi))$-(super)module~$P$ of the form
  $P=\Gamma(\pi_\infty^*(\xi))$ is said to be a \emph{horizontal
  module}.
\end{definition}

\begin{example}[horizontal forms]
  Let $\xi$ be the $q$th exterior degree of the cotangent bundle
  to~$M$. The corresponding horizontal module
  $\Gamma(\pi^*_\infty(\xi))$ is called the module of \emph{horizontal
  forms} and is denoted by $\bar\Lambda^q(\pi)$.  In coordinates,
  horizontal forms are generated by the forms
  $f\,dx_{i_1}\wedge\dots\wedge dx_{i_q}$, $f\in\mathcal{F}(\pi)$.
\end{example}

\begin{definition}
  Let $P$ and $Q$ be $C^\infty(J^\infty(\pi))$-(super)modules.  A map
  $\Delta\colon P\to Q$ is called \emph{$\mathcal{C}$-differential
  operator} (or \emph{horizontal operator}) if it can be written as a
  sum of compositions of $C^\infty(J^\infty(\pi))$-linear maps and
  vector fields of the form~$\hat X$.
\end{definition}

In coordinates, $\mathcal{C}$-differential operators are total
derivatives operators.

\begin{example}[the horizontal de~Rham complex]
  We define the first horizontal de~Rham differential $\bar
  d\colon\mathcal{F}(\pi)\to\bar\Lambda^1(\pi)
  =\Lambda^1(M)\otimes_{C^\infty(M)}\mathcal{F}(\pi)$ by the formula
  $\bar d(f)(X)=\hat X(f)$.  In coordinates, we have $\bar
  d(f)=\sum_iD_i(f)\,dx_i$.
  
  The general horizontal differential $\bar
  d\colon\bar\Lambda^q(\pi)\to\bar\Lambda^{q+1}(\pi)$ is defined by
  the usual rules:
  \begin{align*}
    \bar d\circ\bar d&=0, \\
    \bar d(\omega_1\wedge\omega_2)&=\bar d\omega_1\wedge\omega_2
    +(-1)^q\omega_1\wedge\bar d\omega_2,
    \qquad\omega_1\in\bar\Lambda^q(\pi).
  \end{align*}
  
  The differential $\bar d$ is a $\mathcal{C}$-differential operator.
  
  The cohomology of the horizontal de~Rham complex
  \begin{equation*}
    0 \xrightarrow{} \mathcal{F}(\pi) \xrightarrow{\bar d}
    \bar\Lambda^1(\pi) \xrightarrow{\bar d} \cdots \xrightarrow{\bar d}
    \bar\Lambda^n(\pi)\xrightarrow{}0
  \end{equation*}
  are called \emph{horizontal cohomology} and denoted by $\bar
  H^q(\pi)$.  From Vinogradov's `one line
  theorem'~\cite{Vinogradov:AlGFLFT,Vinogradov:SSAsNDEqAlGFLFTC,Vinogradov:SSLFCLLTNT}
  it follows that $\bar H^q(\pi)=H^q(M)$ for $q\leq n-1$.
\end{example}

All $\mathcal{C}$-differential operators from~$P$ to~$Q$ form a
$C^\infty(J^\infty(\pi))$-(super)module denoted by $\CDiff(P,Q)$.

Clearly, if $P$ and $Q$ are horizontal, then so is $\CDiff(P,Q)$.

Given a horizontal module~$P$, let us define the \emph{horizontal
infinite jet bundle} $\pi_P\colon\bar J^\infty(P)\to J^\infty(\pi)$ as
follows.  If $P$ is even, then the fiber of $\bar J^\infty(P)$ over
$\theta\in J^\infty(\pi)$ consists of equivalence classes, denoted by
$\bar\jmath(p)(\theta)$, of elements~$p\in P$.  Two elements $p_1$ and
$p_2$ are equivalent if their total derivatives of all orders coincide
at~$\theta$.  For a general horizontal supermodule~$P$, we as always
define $\pi_P^0=\pi_{P^0}$ and $\pi_P^1=\pi_{(P^1)^\Pi}^\Pi$.
Correspondingly,
$\Gamma(\pi_P)=\Gamma(\pi_{P^0})\oplus\Gamma(\pi_{(P^1)^\Pi})^\Pi$.

Clearly, the horizontal jet bundle $\pi_P\colon\bar J^\infty(P)\to
J^\infty(\pi)$, $P=\Gamma(\pi_\infty^*(\xi))$, is isomorphic to the
pullback $\pi_\infty^*(\xi_\infty)\colon J^\infty(\xi)\times_M
J^\infty(\pi)\to J^\infty(\pi)$ and, thus, $\Gamma(\pi_P)$ is a
horizontal module.

Similarly to Remark~\ref{sec:cartan-distribution-1}, we have the
natural operator $\bar\jmath_\infty\colon P\to \bar J^\infty(P)$.

For every $\mathcal{C}$-differential operator $\Delta\colon P\to Q$
there exists a unique homomorphism of
$C^\infty(J^\infty(\pi))$-supermodules $h_\Delta\colon\bar
J^\infty(P)\to\bar J^\infty(Q)$ such that the diagram
\begin{equation*}
  \xymatrix{
  P \ar[r]^{\Delta} \ar[d]_{\bar\jmath_\infty} &
  Q \ar[d]^{\bar\jmath_\infty} \\
  \bar J^\infty(P) \ar[r]^{h_\Delta} & \bar J^\infty(Q)
  }
\end{equation*}
is commutative.

Let us recall the definition of adjoint operator.  Consider
$\Delta\in\CDiff(P_1,P_2)$.  The \emph{adjoint operator}
$\Delta^*\in\CDiff(\hat P_2,\hat P_1)$, $\hat
P=\Hom_{\mathcal{F}(\pi)}(P,\bar\Lambda^n(\pi))$, is uniquely defined
by the equality\footnote{Here and below, symbols used at the exponents
of $(-1)$ stand for the corresponding parity.}
\begin{equation}\label{eq:10}
  \langle\hat p,\Delta(p)\rangle
  =(-1)^{\Delta\hat p}\langle\Delta^*(\hat p),p\rangle,
  \qquad\hat p\in\hat P_2,\quad p\in P_1,
\end{equation}
where $\langle\cdot\,,\cdot\rangle$ is the natural pairing $\hat
P\times P\to\bar H^n(\pi)$.

In coordinates, we have
\begin{equation*}
  \matr[\bigg]{\sum_\tau a_{ij}^\tau D_\tau}^*
  =\matr[\bigg]{\sum_\tau(-1)^{\abs{\tau}}D_\tau\circ
  a_{ij}^\tau}^{\text{st}},
\end{equation*}
where $a_{ij}^\tau\in\mathcal{F}(\pi)$, the superscript `st' denotes
the supertransposition, and $D_\tau=D_{i_1}\circ\dots\circ
D_{i_{\abs{\tau}}}$ for $\tau=i_1\dots i_{\abs{\tau}}$.

Equivalently, adjoint operator can be defined using the following
fact.  Consider a horizontal module~$P$ and the natural complex
\begin{multline*}
  0\xrightarrow{}\CDiff(P,\mathcal{F}(\pi))\xrightarrow{}
  \CDiff(P,\bar\Lambda^1(\pi))
  \xrightarrow{}\CDiff(P,\bar\Lambda^2(\pi))\xrightarrow{}\cdots \\
  \xrightarrow{}\CDiff(P,\bar\Lambda^{n-1}(\pi))
  \xrightarrow{}\CDiff(P,\bar\Lambda^n(\pi))\xrightarrow{}0
\end{multline*}
with the differential $\Delta\mapsto\bar d\circ\Delta$.  Denote its
cohomology by $H^q(P)$.  We have
\begin{equation}
  \label{eq:13}
  \begin{cases}
    H^q(P)=0 & \text{ for $0<q<n$} \\
    H^n(P)=\hat P.
  \end{cases}
\end{equation}
Each $\mathcal{C}$-differential operator $\Delta\colon P\to Q$ gives
rise to a cochain map between two such complexes.  The corresponding
map of the $n$th cohomology $\Delta^*\colon\hat Q\to \hat P$ is the
adjoint operator.  Note that the natural projection
$\mu\colon\CDiff(P,\bar\Lambda^n(\pi))\to\hat{P}$ has the form
$\Delta\mapsto\Delta^*(1)$.

Recall the most important properties of adjoint operators:
\begin{enumerate}
\item $\Delta$ and $\Delta^*$ are of equal parity;
\item $\vphantom{\hat{\hat{P}}}(\Delta_1\circ\Delta_2)^*
  =(-1)^{\Delta_1\Delta_2}\Delta_2^*\circ\Delta_1^*$;
\item $\Delta^{**}=\Delta$ (here we identify $\hat{\hat{P}}$ and $P$).
\end{enumerate}

A vector field~$Z$ on $J^\infty(\pi)$ is called \emph{vertical} if
$\eval{Z}_{C^\infty(M)}=0$.  For a horizontal module~$P$ a vertical
field~$Z$ generates a natural action $Z\colon P\to P$, which in
coordinates is the component-wise action.

A vertical vector field~$Z$ is said to be \emph{evolutionary} if
$[Z,\hat X]=0$ for all vector fields $X$ on~$M$.

It is easy to see that evolutionary fields are uniquely determined by
their restrictions to $\mathcal{F}_{\text{\textup{lin}}}(E)$, where
$E$ is the space of the bundle $\pi\colon E\to M$.  Moreover, the map
$Z\mapsto\eval{Z}_{\mathcal{F}_{\text{\textup{lin}}}(E)}$ is a
bijection between the set of all evolutionary fields and
$\Hom_{C^\infty(M)}(\mathcal{F}_{\text{\textup{lin}}}(E),\mathcal{F}(\pi))$.
We identify
$\Hom_{C^\infty(M)}(\mathcal{F}_{\text{\textup{lin}}}(E),\mathcal{F}(\pi))$
with the horizontal module $\Gamma(\pi_\infty^*(\pi))$ and denote it
by~$\kappa(\pi)$.

In coordinate language, the evolutionary field that corresponds to a
vector function $\phi=(\phi^1,\dots,\phi^m)$ has the form
\begin{equation*}
  \re_\phi=\sum_{j,\tau}D_\tau(\phi^j)\pd{}{u^j_\tau}.
\end{equation*}

Let $P$ be a horizontal module.  The \emph{linearization} of an
element $F\in P$ is a $\mathcal{C}$-differential operator
$\ell_F\colon\kappa(\pi)\to P$ defined by the formula
\begin{equation*}
  \ell_F(\phi)=(-1)^{F\phi}\re_\phi(F).
\end{equation*}

Denote by the square brackets the horizontal cohomology class of a
horizontal form.  Since evolutionary fields commute with the
horizontal differential, the cohomology class $[\re_\phi(\omega)]$ for
$\omega\in\bar\Lambda^n(\pi)$ is well defined by~$[\omega]$; denote it
by $\re_\phi([\omega])$.  By~\eqref{eq:10} we have
\begin{equation*}
  \re_\phi([\omega])=[\re_\phi(\omega)]=(-1)^{\phi\omega}[\ell_\omega(\phi)]
  =\langle\phi,\ell_\omega^*(1)\rangle=\langle\phi,\Eu(\omega)\rangle,
\end{equation*}
where $\Eu\colon\bar\Lambda^n(\pi)\to\hat\kappa(\pi)$,
$\Eu(\omega)=\ell_\omega^*(1)$, is the \emph{Euler operator}, which
takes Lagrangians to the corresponding Euler-Lagrange equations.  Of
course, the value $\Eu(\omega)$ is completely determined by the
cohomology class~$[\omega]$.

In coordinates, $\Eu(L\,dx^1\wedge\dots\wedge
dx^n)=(\fd{L}{u^1},\dots,\fd{L}{u^m})$, where
$\fd{L}{u^j}=\sum_\tau(-1)^{\abs{\tau}}D_\tau(\pd{L}{u^j_\tau})$.

\begin{remark}
  \label{sec:jet-bundle-setting-1-line}
  From Vinogradov's `one-line
  theorem'~\cite{Vinogradov:AlGFLFT,Vinogradov:SSAsNDEqAlGFLFTC,%
Vinogradov:SSLFCLLTNT} it follows that
  \begin{enumerate}
  \item $\ker\Eu/\bar d(\bar\Lambda^{n-1}(\pi))=H^n(M)$;
  \item $\psi\in\im\Eu$ if and only if $\ell_\psi^*=\ell_\psi$.
  \end{enumerate}
\end{remark}

\subsection{Differential equations}
\label{sec:diff-equat}

Again, consider an element $F$ of a horizontal module~$P$.  The locus
\begin{equation*}
  \mathcal{E}^\infty=\{\,\bar\jmath_\infty(F)=0\,\}\subset J^\infty(\pi)
\end{equation*}
is called \emph{differential equation} defined by~$F$.  We assume that
the natural map $\mathcal{E}^\infty\to M$ is a subbundle of the bundle
$\pi\colon J^\infty(\pi)\to M$.  The restriction of the Cartan
distribution to~$\mathcal{E}^\infty$ is denoted by
$\mathcal{C}(\mathcal{E})$.  Clearly,
$\dim\mathcal{C}(\mathcal{E})=\dim\mathcal{C}=n$.

\begin{example}[evolution equations]
  Consider the bundle $\bar\pi\colon E\times\mathbb{R}\to
  M\times\mathbb{R}$.  Denote the coordinate along~$\mathbb{R}$
  by~$t$.  Then $\re_\Phi=D_t-\pd{}{t}$ is a canonical evolutionary
  field on~$J^\infty(\bar\pi)$.  In coordinates,
  $\Phi=(u^1_t,\dots,u^m_t)$.  Let $\re_{\phi(t)}$ be a family of
  evolutionary fields on~$J^\infty(\pi)$.  The equation
  $\mathcal{E}^\infty\subset J^\infty(\bar\pi)$ given by the element
  $F=\Phi-\phi(t)\in\kappa(\bar\pi)$ is called \emph{evolution
  equation}.  In coordinates, it has the form $u_t=\phi(t)$.  Note
  that $\mathcal{E}^\infty=J^\infty(\pi)\times\mathbb{R}$, with the
  Cartan distribution generated by that on $J^\infty(\pi)$ and the
  field $D_t=\pd{}{t}+\re_{\phi(t)}$.
\end{example}

The restriction of the linearization $\ell_F$ to the equation
$\mathcal{E}^\infty$ is called the \emph{linearization
of~$\mathcal{E}^\infty$} and is denoted by
$\ell_\mathcal{E}\colon\kappa\to P$, where
$\kappa=\eval{\kappa(\pi)}_{\mathcal{E}^\infty}$.

An evolutionary field tangent to $\mathcal{E}^\infty$ is said to be a
\emph{symmetry} of the equation.  Obviously, $\re_\phi$ is a symmetry
if and only if $\ell_\mathcal{E}(\phi)=0$, $\phi\in\kappa$.

The horizontal de~Rham complex on $J^\infty(\pi)$ can be restricted
to~$\mathcal{E}^\infty$.  Its cohomology are called \emph{horizontal
cohomology of equation~$\mathcal{E}^\infty$} and denoted by~$\bar
H^q(\mathcal{E})$.  Elements of $\bar H^{n-1}(\mathcal{E})/H^{n-1}(M)$
are \emph{conservation laws} of~$\mathcal{E}^\infty$.  If the equation
at hand satisfies the conditions of Vinogradov's `two-line
theorem'~\cite{Vinogradov:AlGFLFT,Vinogradov:SSAsNDEqAlGFLFTC,%
Vinogradov:SSLFCLLTNT} then there is an inclusion $i\colon\bar
H^{n-1}(\mathcal{E})/H^{n-1}(M)\to\ker\ell^*_\mathcal{E}$. The element
$i(\eta)\in\ker\ell^*_\mathcal{E}\subset\hat P$ that corresponds to
conservation law~$\eta$ is called its \emph{generating function}.

In particular, evolution equations satisfy the conditions of the
two-line theorem.  In this case, $i(\eta)=\Eu(\eta_0)$, where
$\eta=\eta_0+\eta_1\wedge dt$,
$\eta_0\in\bar\Lambda^{n-1}(\mathcal{E})$,
$\eta_1\in\bar\Lambda^{n-2}(\mathcal{E})$.  Thus, to find conservation
laws of an evolution equation one has to solve the equation
$\ell_\mathcal{E}^*(\psi)=0$ and choose those solutions $\psi$ that
fulfill the condition $\ell_\psi=\ell^*_\psi$.

Let $\mathcal{E}^\infty$ and $\tilde{\mathcal{E}}^\infty$ be two
differential equations. A surjective map $\tau\colon
\tilde{\mathcal{E}}\to\mathcal{E}^\infty$ is called \emph{covering} if
it preserves the Cartan distribution.

\begin{example}
  A horizontal jet bundle $\pi_P\colon\bar J^\infty(P)\to
  J^\infty(\pi)$ is a covering.  A generalization of this example will
  be discussed in Section~\ref{sec:delta-coverings}.
\end{example}

\begin{example}\label{ex:cover}
  Let $\mathcal{E}^\infty$ be given by an element $F$. Consider
  equation $\tilde{\mathcal{E}}^\infty$
  \begin{equation}
    \label{eq:14}
    \begin{cases}
      F=0, \\
      \pd{r}{x_1}=g_1, \\
      \newlength{\DotsAs}\settowidth{\DotsAs}{$\pd{r}{x_1}=g_1,$}
      \hbox to\DotsAs{\dotfill} \\
      \pd{r}{x_n}=g_n,
    \end{cases}
  \end{equation}
  where $g_1,\dots,g_n$ are functions of $x_1,\dots,x_n$ and
  $u_\sigma^j$.  If the compatibility condition
  \begin{equation}
    \label{eq:12}
    \text{$D_i(g_j)=D_j(g_i)$ on~$\mathcal{E}^\infty$}\quad\iff\quad
    \text{$\eval{\bar d\eta}_{\mathcal{E}^\infty}=0$\qquad
    $\eta=\textstyle\sum_i g_i\,dx_i\in\bar\Lambda^1(\mathcal{E})$}
  \end{equation}
  holds true, then the natural projection
  $\tau\colon\tilde{\mathcal{E}}^\infty\to\mathcal{E}^\infty$ is
  epimorphic.  Obviously, $\tau$ preserves the Cartan distribution, so
  that it is a covering.  Thus, each closed horizontal $1$-form gives
  rise to a covering of the form~\eqref{eq:14}.  In particular, when
  $n=2$, condition~\eqref{eq:12} means that $\eta$ represents a
  conservation law of~$\mathcal{E}^\infty$.  The new dependent
  variable~$r$ is called \emph{nonlocal variable}.
  
  In a similar way, we can define a covering
  over~$\tilde{\mathcal{E}}^\infty$ corresponding to a closed one-form
  on this equation, etc.  In this manner we construct particular
  coverings in Sections~\ref{sec:korteweg-de-vries}--\ref{sec:appl}.
\end{example}

Clearly, each $\mathcal{C}$-differential operator~$\Delta$
on~$\mathcal{E}^\infty$ can be lifted to a $\mathcal{C}$-differential
operator~$\tilde\Delta$ on~$\tilde{\mathcal{E}}^\infty$.  In
particular, we have the operator $\tilde\ell_{\mathcal{E}}$ on
$\tilde{\mathcal{E}}^\infty$. A symmetry of
$\tilde{\mathcal{E}}^\infty$ is called a \emph{nonlocal symmetry} of
$\mathcal{E}^\infty$ in the covering under consideration. Solutions of
the equation $\tilde\ell_{\mathcal{E}}(\phi)=0$ are called
\emph{shadows of nonlocal symmetries} of~$\mathcal{E}^\infty$ in this
covering. In a similar way, since the horizontal de~Rham differential
is a $\mathcal{C}$-differential operator, we can lift the horizontal
de~Rham complex to $\tilde{\mathcal{E}}^\infty$ and construct the
theory of \emph{nonlocal conservation laws} in our covering.
Solutions of the equation $\tilde{\ell_\mathcal{E}^*}(\psi)=0$ are
called \emph{nonlocal generating functions}.

\section{Variational Schouten bracket}
\label{sec:vari-scho-brack}

We start with a super version of Kupershmidt's \emph{cotangent bundle
to a vector bundle}~\cite{Kupershmidt:GJBSLHF}.

For a vector bundle~$\pi\colon E\to M$, $\dim M=n$, we consider the
bundle $\hat\pi\colon \hat E=E^*\otimes_M\Lambda^n(T^*M)\to M$, where
$E^*\to M$ is the dual bundle to~$E\to M$, and the superbundle
$\mathcal{K}$: $\mathcal{K}^0=\pi$ (even subbundle),
$\mathcal{K}^1=\hat\pi$ (odd subbundle).

The superbundle $\mathcal{K}_\infty\colon J^\infty(\mathcal{K})\to M$
\begin{equation*}
  \xymatrix {
  J^\infty(\mathcal{K}) \ar[dr] ^{\pi_\infty^*(\hat{\pi}_\infty)}
  \ar[dd]_{\mathcal{K}_\infty} \\
  & J^\infty(\pi) \ar[dl]^{\pi_\infty} \\
  M
  }
\end{equation*}
is called the \emph{cotangent bundle} of the bundle~$\pi$.  It is
clear that
\begin{equation*}
  J^\infty(\mathcal{K})
  =\bar J^\infty(\Gamma(\pi_\infty^*(\hat{\pi}_\infty))).
\end{equation*}

Denote by $p^j$, $j=1,\dots,m$, the fiber coordinates in $\hat{E}$
dual to~$u^j$ with respect to a volume form on~$M$ (they are sometimes
called `antifields'). Then $x_i$, $u^j_\tau$, $p^j_\tau$ will be the
coordinates in $J^\infty(\mathcal{K})$, with $x_i$, $u^j_\tau$ being
even and $p^j_\tau$ being odd.

It is clear that $\kappa(\mathcal{K})^0=\kappa_\mathcal{K}$ and
$\kappa(\mathcal{K})^1=\hat\kappa_\mathcal{K}^\Pi$, where
$\kappa_\mathcal{K} = \Gamma(\mathcal{K}^*(\pi))$.

Define the \emph{variational Schouten bracket
\textup{(}antibracket\textup{)}} on the space~$\bar H^n(\mathcal{K})$
by putting
\begin{equation}\label{eq:pb}
  [\![F,H]\!] = \langle\Eu(H),\alpha(\Eu(F))\rangle,
  \qquad F,\ H\in\bar H^n(\mathcal{K}),
\end{equation}
where $\Eu$ is the Euler operator and the
operator~$\alpha\colon\hat{\kappa}(\mathcal{K}) \to
\kappa(\mathcal{K})$ acts according to the formula
$\alpha(\psi,\phi)=(\phi,-\psi)$ for $\phi\in\kappa_\mathcal{K}$ and
$\psi\in\hat\kappa_\mathcal{K}$. In coordinates, we get
\begin{equation*}
  [\![F,H]\!]=\sum_j\Bigl[\fd{H}{u^j}\fd{F}{p^j}
  -(-1)^{(F+1)(H+1)}\fd{F}{u^j}\fd{H}{p^j}\Bigr].
\end{equation*}

It is readily seen that the variational Schouten bracket defines a Lie
superalgebra structure on~$\bar H^n(\mathcal{K})$:
\begin{gather*}
  [\![H,F]\!]=-(-1)^{(F+1)(H+1)}[\![F,H]\!], \\[1ex]
  \begin{split}
    (-1)^{(F+1)(H+1)}[\![[\![F,G]\!],H]\!]
    &+(-1)^{(F+1)(G+1)}[\![[\![G,H]\!],F]\!] \\
    &+(-1)^{(G+1)(H+1)}[\![[\![H,F]\!],G]\!]=0.
  \end{split}
\end{gather*}

\begin{remark}
  A different concept of the Schouten bracket (acting on a different
  space) the reader can find in~\cite[p.~226]{Vinogradov:CAnPDEqSC}.
\end{remark}

Denote by $\CDiffskew{k}(P,Q)$ the module of $k$-linear skew-symmetric
$\mathcal{C}$-differential operators $P\times\dots\times P\to Q$. The
subset $\CDiffskewad{k}(P,\hat P)\subset\CDiffskew{k}(P,\hat P)$
consists of skew-adjoint in each argument operators.

Let us define multiplication
\begin{equation*}
  \CDiffskew{k}(P,\mathcal{F}(\pi))
  \times\CDiffskew{l}(P,\mathcal{F}(\pi))
  \xrightarrow{}\CDiffskew{k+l}(P,\mathcal{F}(\pi))
\end{equation*}
by setting
\begin{equation*}
  (\Delta_1\Delta_2) (p_1,\dots ,p_{k+l}) \\
  =\sum_{\sigma\in S_{k+l}^k}(-1)^{\sigma}
  \Delta_1(p_{\sigma(1,k)})\Delta_2(p_{\sigma(k+1,k+l)}),
\end{equation*}
where $\Delta_1\in\CDiffskew{k}(P,\mathcal{F}(\pi))$,
$\Delta_2\in\CDiffskew{l}(P,\mathcal{F}(\pi))$, $S_n^i\subset S_n$ is
the set of all $(i,n-i)$-unshuffles (\cite{LadaStasheff:InLAlP}),
i.e., all permutations~$\sigma\in S_n$ such that
$\sigma(1)<\sigma(2)<\dots<\sigma(i)$ and
$\sigma(i+1)<\sigma(i+2)<\dots<\sigma(n)$, $(-1)^\sigma$ is the sign
of permutation~$\sigma$, and $p_{\sigma(k_1,k_2)}$ stands for
$p_{\sigma(k_1)},\dots,p_{\sigma(k_2)}$.

Next, since by definition elements of $\mathcal{F}(\hat{\pi})$ are
identified with differential operators from $\Gamma(\hat{\pi})$ to
$C^\infty(M)$, we have the natural inclusion
$\CDiff(\hat{\kappa}(\pi),\mathcal{F}(\pi)) \to
\mathcal{F}(\mathcal{K})$, which uniquely prolongs to the isomorphism
of algebras
\begin{equation*}
  \CDiffskew{*}(\hat{\kappa}(\pi),\mathcal{F}(\pi))
  \to\mathcal{F}(\mathcal{K}).
\end{equation*}

Using~\eqref{eq:13} we can show in a standard way that
\begin{equation}
  \label{eq:9}
  \bar H^n(\mathcal{K})=\CDiffskewad{*}(\hat\kappa(\pi),\kappa(\pi))
  \oplus\bar{H}^n(\pi).
\end{equation}

Below, we use the shorthand notation $\kappa=\kappa(\pi)$.

Now, following~\cite{IgoninVerbovetskyVitolo:FLVDOp}, we want to
compute the variational Schouten bracket in terms of skew-adjoint
$\mathcal{C}$-differential operators.

To this end, note that from the definition of the Euler operator it
follows that its restriction
\begin{equation*}
  \eval{\Eu}_{\CDiffskew{k}(\hat{\kappa},\bar\Lambda^n(\pi))}
  \colon\CDiffskew{k}(\hat{\kappa},\bar\Lambda^n(\pi))\to
  \CDiffskew{k}(\hat{\kappa},\hat{\kappa})\oplus
  \CDiffskew{k-1}(\hat{\kappa},\kappa),
\end{equation*}
has the form
$\eval{\Eu}_{\CDiffskew{k}(\hat{\kappa},\bar\Lambda^n(\pi))}(\Delta)
=(\eta(\Delta),(-1)^{k-1}\mu(\Delta))$, where
\begin{align*}
  \eta(\Delta)(\psi_1,\dots,\psi_k)
  &=\ell^*_{\Delta,\psi_1,\dots,\psi_k}(1), \\
  \mu(\Delta)(\psi_1,\dots,\psi_{k-1})
  &=(\Delta(\psi_1,\dots,\psi_{k-1}))^*(1),
\end{align*}
$\psi_i\in\hat\kappa$,
$\ell_{\Delta,\psi_1,\dots,\psi_k}(\phi)=\re_\phi(\Delta)(\psi_1,\dots,\psi_k)$.

In coordinates, $\eta=(\fd{}{u^1},\dots ,\fd{}{u^m})$ and
$\mu=(-1)^{k-1}(\fd{}{p^1},\dots ,\fd{}{p^m})$.

We can rewrite~$\eta$ in the following form:

\begin{equation*}
  \eta(\Delta)=\tilde\eta(\mu(\Delta)),\qquad
  \Delta\in\CDiffskew{k}(\hat{\kappa},\bar\Lambda^n(\pi)),
\end{equation*}
where
\begin{equation*}
  \tilde\eta(\square)(\psi_1,\dots,\psi_k)
  =\ell^*_{\square,\psi_1,\dots,\psi_{k-1}}(\psi_k),
  \quad\square\in\CDiffskewad{k-1}(\hat\kappa,\kappa).
\end{equation*}

Indeed, take the equality
\begin{equation*}
  [\Delta(\psi_1,\dots,\psi_k)]
  =\langle\square(\psi_1,\dots,\psi_{k-1}),\psi_k\rangle,\qquad
  \square=\mu(\Delta),
\end{equation*}
and apply $\re_\phi$ to both sides. This yields
\begin{equation*}
  [\re_\phi(\Delta)(\psi_1,\dots,\psi_k)]
  =\langle\re_\phi(\square)(\psi_1,\dots,\psi_{k-1}),\psi_k\rangle,
\end{equation*}
and so
\begin{equation*}
  \langle\phi,\ell^*_{\Delta,\psi_1,\dots,\psi_k}(1)\rangle
  =\langle\phi,\ell^*_{\square,\psi_1,\dots,\psi_{k-1}}(\psi_k)\rangle.
\end{equation*}

Thus, for the variational Schouten bracket of two operators
$A\in\CDiffskewad{k-1}(\hat\kappa,\kappa)$ and
$B\in\CDiffskewad{l-1}(\hat\kappa,\kappa)$ we have
\begin{multline*}
  \langle[\![A,B]\!](\psi_1,\dots,\psi_{k+l-2}),\psi_{k+l-1}\rangle \\
  =[((-1)^{k-1}\tilde\eta(B)A-(-1)^{k(l-1)}\tilde\eta(A)B)
  (\psi_1,\dots,\psi_{k+l-1})] \\
  =(-1)^{k-1}\sum_{\sigma\in S_{k+l-1}^l}(-1)^\sigma
  \langle\ell^*_{B,\psi_{\sigma(1,l-1)}}
  (\psi_{\sigma(l)}),A(\psi_{\sigma(l+1,k+l-1)})\rangle \\
  -(-1)^{k(l-1)}\sum_{\sigma\in S_{k+l-1}^k}(-1)^\sigma
  \langle\ell^*_{A,\psi_{\sigma(1,k-1)}}
  (\psi_{\sigma(k)}),B(\psi_{\sigma(k+1,k+l-1)})\rangle.
\end{multline*}

Here and below we assume that $S_n^i=\varnothing$ if $i<0$ or $i>n$.

Let us split the sums obtained into two parts depending on whether
$\sigma(k+l-1)=k+l-1$ or not:
\begin{multline*}
  \langle[\![ A,B ]\!](\psi_1,\dots,
  \psi_{k+l-2}),\psi_{k+l-1}\rangle \\
  =(-1)^{k-1}\sum_{\sigma\in S_{k+l-2}^l}(-1)^\sigma
  \langle\ell^*_{B,\psi_{\sigma(1,l-1)}}
  (\psi_{\sigma(l)}),A(\psi_{\sigma(l+1,k+l-2)},\psi_{k+l-1})\rangle \\
  +\sum_{\sigma\in S_{k+l-2}^{l-1}}(-1)^\sigma
  \langle\ell^*_{B,\psi_{\sigma(1,l-1)}}
  (\psi_{k+l-1}),A(\psi_{\sigma(l,k+l-2)})\rangle \\
  -(-1)^{k(l-1)}\sum_{\sigma\in S_{k+l-2}^k}(-1)^\sigma
  \langle\ell^*_{A,\psi_{\sigma(1,k-1)}}
  (\psi_{\sigma(k)}),B(\psi_{\sigma(k+1,k+l-2)},\psi_{k+l-1})\rangle
  \\
  -(-1)^{(k-1)(l-1)}\sum_{\sigma\in S_{k+l-2}^{k-1}}(-1)^\sigma
  \langle\ell^*_{A,\psi_{\sigma(1,k-1)}}
  (\psi_{k+l-1}),B(\psi_{\sigma(k,k+l-2)})\rangle.
\end{multline*}

Thus, we have
\begin{multline}\label{eq:1}
  [\![A,B]\!](\psi_1,\dots,\psi_{k+l-2}) =\sum_{\sigma\in
  S_{k+l-2}^{l-1}}(-1)^\sigma
  \ell_{B,\psi_{\sigma(1,l-1)}}(A(\psi_{\sigma(l,k+l-2)})) \\
  -(-1)^{(k-1)(l-1)}\sum_{\sigma\in S_{k+l-2}^k}(-1)^\sigma
  B(\ell^*_{A,\psi_{\sigma(1,k-1)}}(\psi_{\sigma(k)}),
  \psi_{\sigma(k+1,k+l-2)}) \\
  -(-1)^{(k-1)(l-1)}\sum_{\sigma\in S_{k+l-2}^{k-1}}(-1)^\sigma
  \ell_{A,\psi_{\sigma(1,k-1)}}(B(\psi_{\sigma(k,k+l-2)})) \\
  +\sum_{\sigma\in S_{k+l-2}^l}(-1)^\sigma
  A(\ell^*_{B,\psi_{\sigma(1,l-1)}}
  (\psi_{\sigma(l)}),\psi_{\sigma(l+1,k+l-2)}).
\end{multline}

From the definition it immediately follows that
\begin{equation*}
  [\![A,\omega]\!](\psi_1,\dots,\psi_{k-2})=A(\Eu(\omega),\psi_1,\dots,\psi_{k-2})
\end{equation*}
for $\omega\in\bar H^n(\pi)$; in particular,
$[\![\phi,\omega]\!]=\langle\phi,\Eu(\omega)\rangle
=\re_\phi(\omega)$.

\section{Hamiltonian evolution equations}
\label{sec:hamilt-evol-equat}

An operator $A\in\CDiff(\hat\kappa,\kappa)$ is called
\emph{Hamiltonian} if $[\![A,A]\!]=0$. As in the classical Hamiltonian
formalism, a Hamiltonian operator defines a Lie algebra structure
on~$\bar H^n(\pi)$ via the \emph{Poisson bracket}
\begin{equation*}
  \{\omega_1,\omega_2\}_A
  =\langle A(\Eu(\omega_1)),\Eu(\omega_2)\rangle.
\end{equation*}

\begin{remark}
  Hamiltonian operators are uniquely determined by the corresponding
  Poisson brackets.
\end{remark}

\begin{remark}
  A Hamiltonian operator~$A$ gives rise to a complex
  \begin{equation}
    \label{eq:11}
    0\xrightarrow{}\bar H^n(\pi)\xrightarrow{\partial_A}
    \kappa\xrightarrow{\partial_A}\CDiff(\hat\kappa,\kappa)
    \xrightarrow{\partial_A}\CDiffskewad{2}(\hat\kappa,\kappa)
    \xrightarrow{\partial_A}\cdots,
  \end{equation}
  where $\partial_A(\Delta)=[\![A,\Delta]\!]$, called the
  \emph{Hamiltonian complex}.
\end{remark}

Formula~\eqref{eq:1} yields a well-known criterion for checking a
skew-adjoint operator to be Hamiltonian (see,
e.g.,~\cite{KrasilshchikVinogradov:SCLDEqMP,KrasilshchikVerbovetsky:HMEqMP}):
\begin{equation*}
  [\![A,A]\!](\psi_1,\psi_2)=-\ell_{A,\psi_1}(A(\psi_2))
  +\ell_{A,\psi_2}(A(\psi_1))-A(\ell^*_{A,\psi_1}(\psi_2))=0.
\end{equation*}

Another practical way to check the Hamiltonian property of an operator
is to use formula~\eqref{eq:pb}. In coordinates, it gives:
\begin{align}
  \label{eq:2}
  &\sum_j\fd{W_A}{u^j}\fd{W_A}{p^j}=0\ \text{modulo total derivatives} \\
  \text{or}\ 
  \Eu\Bigl(&\sum_j\fd{W_A}{u^j}\fd{W_A}{p^j}\Bigr)=0,\nonumber
\end{align}
where $W_A\in\bar H^n(\mathcal{K})$ is the element that corresponds to
the operator $2A$ under the isomorphism~\eqref{eq:9}. In coordinates,
the element~$W_A$ for an operator $\sum_\tau a^{ij}_\tau D_\tau$ has
the form $W_A=\sum_{i,j,\tau} a^{ij}_\tau p^j_\tau p^i$.

The condition for two Hamiltonian operators $A$ and $B$ to be a
Hamiltonian pair, i.e., $[\![A,B]\!]=0$, is
\begin{equation*}
  \sum_j\Eu\Bigl(\fd{W_A}{u^j}\fd{W_B}{p^j}
  +\fd{W_B}{u^j}\fd{W_A}{p^j}\Bigr)=0.
\end{equation*}

Note that the skew-adjointness in terms of~$W_A$ amounts to the
equality
\begin{equation}
  \label{eq:3}
  \sum_j\fd{W_A}{p^j}p^j=-2W_A.
\end{equation}

Let $A$ be a Hamiltonian operator.  Evolution equation $u_t=f$ is said
to be \emph{Hamiltonian} with respect to~$A$ if
\begin{equation}
  \label{eq:4}
  A_t-[\![A,f]\!]=0,
\end{equation}
where $A_t=\pd{A}{t}$ (both $A$ and $f$ can depend on the
parameter~$t$).

If $A$ does not depend on~$t$, then for each $H\in\bar H^n(\pi)$ the
evolution equation $u_t=A(\Eu(H))$ is a Hamiltonian evolution
equation.  The element $H\in\bar H^n(\pi)$ is called the
\emph{Hamiltonian}.  Notice that in this case condition~\eqref{eq:4}
means that $f$ is a $1$-cocycle in the Hamiltonian
complex~\eqref{eq:11}.  A Hamiltonian~$H$ exists if and only if $f$ is
a coboundary.

If a Hamiltonian~$H$ exists and does not depend on~$t$, then we have
\begin{equation*}
  D_t(H)=\re_f(H)=\re_{A(\Eu(H))}(H)=\{H,H\}_A=0.
\end{equation*}
Thus, there exists a conservation law given by $\eta_0+\eta_1\wedge
dt$, $\eta_0\in\bar\Lambda^n(\pi)$, $\eta_1\in\bar\Lambda^{n-1}(\pi)$,
such that $[\eta_0]=H\in\bar{H}^n(\pi)$, where $[\eta_0]$ is the
cohomology class of~$\eta_0$ in $\bar H^n(\pi)$.  In other words, the
generating function of this conservation law equals~$\Eu(H)$.  This
conservation law is called the \emph{conservation law of energy}.

\begin{theorem}
  Let $\mathcal{E}^\infty$ be an evolution equation $u_t=f$ which is
  Hamiltonian with respect to a Hamiltonian operator~$A$. Then we have
  \begin{equation}
    \label{eq:5}
    \ell_\mathcal{E}\circ A+A\circ\ell^*_\mathcal{E}=0.
  \end{equation}
\end{theorem}

\begin{proof}
  By~\eqref{eq:1},
  \begin{equation*}
    (A_t-[\![A,f]\!])(\psi)
    =A_t(\psi)+\ell_{A,\psi}(f)-A(\ell^*_f(\psi))-\ell_f(A(\psi)),
  \end{equation*}
  thus
  \begin{equation*}
    A_t-[\![A,f]\!]=A_t+\re_f(A)-A\circ\ell^*_f-\ell_f\circ A.
  \end{equation*}
  Hence, $(\pd{}{t}+\re_f-\ell_f)\circ
  A-A\circ(\pd{}{t}+\re_f+\ell^*_f)=0$.  It remains to note that
  $\ell_\mathcal{E}=D_t-\ell_f=\pd{}{t}+\re_f-\ell_f$.
\end{proof}

\begin{remark}
  For equations possessing a Hamiltonian, relation~\eqref{eq:5} can be
  found elsewhere (see, e.g.~\cite{Olver:ApLGDEq}).
\end{remark}

We call solutions of~\eqref{eq:5} \emph{variational bivectors on the
equation} under consideration; Hamiltonian operators that make a given
equation Hamiltonian are, thus, special variational bivectors on the
equation. Obviously, variational bivectors (and, in particular,
Hamiltonian operators) take generating functions of conservation laws
of the equation at hand to symmetries of this equation.

\begin{proposition}
  Let $\mathcal{E}^\infty$ be an evolution equation $u_t=f$.  If two
  operators $A$\textup{,} $A'\in\CDiff(\hat\kappa,\kappa)$ satisfy the
  equation
  \begin{equation}
    \label{eq:6}
    \ell_\mathcal{E}\circ A+A'\circ\ell^*_\mathcal{E}=0,
  \end{equation}
  then $A'=A$.
\end{proposition}

\begin{proof}
  Rewrite~\eqref{eq:6} in the form
  \begin{equation*}
    (D_t-\ell_f)\circ A-A'\circ(D_t+\ell^*_f)=0.
  \end{equation*}
  Commute the right-hand side of this equality with the operator of
  multiplication by~$t$. This gives $A'=A$.
\end{proof}

\section{$\Delta$-coverings}
\label{sec:delta-coverings}
In this section, we describe a construction, that reduces solution of
equation~\eqref{eq:5} to finding shadows of nonlocal symmetries in a
special covering over~$\mathcal{E}$ (the
$\ell_\mathcal{E}^*$-covering).

Let $\mathcal{E}^\infty$ be a differential equation, and $\Delta\colon
P\to Q$ be a $\mathcal{C}$-differential operator between two
horizontal modules $P$ and $Q$ over~$\mathcal{E}^\infty$. Consider the
homomorphism $h_\Delta\colon\bar J^\infty(P)\to\bar J^\infty(Q)$ that
corresponds to~$\Delta$. If $K_\Delta=\ker h_\Delta\subset\bar
J^\infty(P)$ is a subbundle of~$\bar J^\infty(P)$, then
$k_\Delta=\eval{\pi_P}_{K_\Delta}\colon K_\Delta\to\mathcal{E}^\infty$
is a covering. We call it \emph{$\Delta$-covering}.

In terms of local coordinates, if $\Delta=\matr[\big]{\sum_\tau
a^{ij}_\tau D_\tau}$ and $w^j$ are fiber coordinates of~$\alpha$,
where $\alpha$ is such that $P=\Gamma(\alpha)$, then $\Delta$-covering
is defined by the equations
\begin{equation}
  \label{eq:7}
  \sum_{\tau j}a^{ij}_\tau w^j_\tau=0.
\end{equation}

We can think of fibers of $\Delta$-covering as even or odd. Here we
prefer the latter viewpoint, so that $k_\Delta$ is a superbundle.

$\Delta$-coverings are useful mainly due to the following obvious
fact.

\begin{proposition}
  \label{sec:delta-covering-1}
  Let $R=\Gamma(\gamma)$ be a horizontal module
  over~$\mathcal{E}^\infty$. Then there is an isomorphism
  \begin{equation*}
    \Gamma_{\text{\textup{lin}}}(k^*_\Delta(\gamma))
    =\CDiff(P,R)\big/\{\,V\in\CDiff(P,R)\mid V=\square\circ\Delta,
    \square\in\CDiff(Q,R)\,\},
  \end{equation*}
  where $\Gamma_{\text{\textup{lin}}}$ denotes space of fiberwise
  linear sections.
\end{proposition}
\begin{proof}
  The isomorphism takes $s\in\Gamma_{\text{lin}}(k^*_\Delta(\gamma))$
  to the equivalence class of the operator~$V_s\colon P\to R$ given by
  the formula
  \begin{equation*}
    V_s(p)=\tilde s\circ\bar\jmath_\infty(p),
  \end{equation*}
  where $\tilde s$ is an extension of~$s$ to $\bar J^\infty(P)$, $p\in
  P$.
\end{proof}

In coordinate language, $D_\tau$ at the $j$th component of the
operator goes to~$w^j_\tau$.

Now suppose that we are given a $\mathcal{C}$-differential operator
$\nabla\colon R\to R'$ over~$\mathcal{E}^\infty$. Let us lift it
on~$K_\Delta$ and consider $\ker\nabla$. In view of
Proposition~\ref{sec:delta-covering-1}, we can identify fiberwise
linear elements of $\ker\nabla$ with
solutions~$V\in\CDiff(P,R)\big/\{\,\square\circ\Delta\,\}$ of the
equation
\begin{equation}
  \label{eq:8}
  \nabla\circ V=V'\circ\Delta
\end{equation}

\emph{Thus\textup{,}~\eqref{eq:8} amounts to the equation
$\nabla(s)=0$ on the $\Delta$-covering}.

In particular, equation~\eqref{eq:6} is equivalent to the equation
$\ell_\mathcal{E}(\phi)=0$ on the $\ell^*_\mathcal{E}$-covering, where
$\phi$ is fiberwise linear vector function. Note that in this case the
$\ell^*_\mathcal{E}$-covering can be identified with the cotangent
bundle~$J^\infty(\mathcal{K})\times\mathbb{R}$. Under this
identification, the Cartan planes on~$K_{\ell^*_\mathcal{E}}$ are
spans of the Cartan planes on~$J^\infty(\mathcal{K})$ and
$D_t=\pd{}{t}+\re_{\tilde{f}}$, where $\tilde{f}=(f,\ell_f^*(w))$, if
the equation at hand is $u_t=f$.  Moreover, the equivalence classes of
operators from Proposition~\ref{sec:delta-covering-1} are in
one-to-one correspondence with $\mathcal{C}$-differential operators
on~$J^\infty(\pi)$.
\begin{remark}
  From the above said, we see that Hamiltonian operators are shadows
  of nonlocal symmetries in the $\ell_\mathcal{E}^*$-covering.
\end{remark}
\begin{remark}
  Solutions~$V\in\CDiff(P,R)\big/\{\,\square\circ\Delta\,\}$ of
  equation~\eqref{eq:8} can be found straightforwardly and the
  computations will be essentially the same as when one solves the
  equation $\nabla(s)=0$ on the $\Delta$-covering. Nevertheless, in
  our computations we prefer the second approach, because in the case
  $\nabla=\ell_{\mathcal{E}}$ it reduces the problem to finding
  shadows of nonlocal symmetries (see above) for which efficient
  software exists.
\end{remark}

\section{The Korteweg-de Vries equation}\label{sec:korteweg-de-vries}
Here we show how the above introduced techniques work with a simple
and well-known example of the KdV equation
\begin{equation}\label{eq:KdV}
  u_t=u_{xxx}+uu_x.
\end{equation}

Local coordinates in $\mathcal{E}^\infty$ are
\begin{equation*}
  x,\ t,\ u=u_0,\dots,u_k,\dotsc,
\end{equation*}
where $u_k=\pd{^ku}{x^k}$ (similar notation is used in the subsequent
sections as well). In these coordinates, the total derivatives are
\begin{align*}
  D_x&=\pd{}{x}+\sum_{k\ge0}u_{k+1}\pd{}{u_k},\\
  D_t&=\pd{}{t}+\sum_{k\ge0}D_x^k(u_3+uu_1)\pd{}{u_k}.
\end{align*}

\subsubsection*{The $\ell_{\mathcal{E}}^*$-covering}
The linearization operator for~\eqref{eq:KdV} is
\begin{equation*}
  \ell_{\mathcal{E}}=D_t-D_x^3-uD_x-u_1,
\end{equation*}
while the adjoint is expressed by the formula
\begin{equation*}
  \ell_{\mathcal{E}}^*=-D_t+D_x^3+uD_x.
\end{equation*}

Following the general scheme, we construct the
$\ell_{\mathcal{E}}^*$-covering by introducing the odd variables
$p=p_0$, $p_k=D_x^k(p)$ that satisfy the equation
\begin{equation*}
  p_t=p_3+up_1.
\end{equation*}

\subsubsection*{Solving the defining equation}
Let us now extend the total derivatives up to the total derivatives on
the $\ell^*_\mathcal{E}$-covering
\begin{align*}
  \tilde{D}_x&=D_x+\sum_{k\ge0}p_{k+1}\pd{}{p_k},\\
  \tilde{D}_t&=D_t+\sum_{k\ge0}\tilde{D}_x^k(p_3+up_1)\pd{}{p_k}.
\end{align*}
Then, solving the equation $\ell_\mathcal{E}(F)=0$, that is,
\begin{equation}
  \label{eq:KdV-ham}
  \tilde{D}_t(F)=\tilde{D}_x^3(F)+u\tilde{D}_x(F)+u_1F,
\end{equation}
with respect to the function $F=\sum_iF_ip_i$, where
$F_i=F_i(x,t,u,\dots,u_k)$, we obtain two independent solutions
\begin{equation*}
  F^0=p_1,\qquad F^1=p_3+\frac{2}{3}up_1+\frac{1}{3}u_1p_0
\end{equation*}
to which there correspond two $\mathcal{C}$-differential operators
\begin{equation*}
  A^0=D_x,\qquad A^1=D_x^3+\frac{2}{3}uD_x+\frac{1}{3}u_1,
\end{equation*}
the classical Hamiltonian structures for the KdV equation.

\subsubsection*{The Hamiltonianity test}
To demonstrate how the method works, we shall check the Hamiltonianity
of the operators $A^0$ and $A^1$ in a straightforward way. Obviously,
both operators are skew-adjoint.

For $A^0$, the corresponding bivector $W_0=W_{A^0}$ is
\begin{equation*}
  W_0=F^0p_0=p_1p_0.
\end{equation*}
Since $\fd{W_0}{u}=0$, we get
\begin{equation*}
  \Eu\Bigl(\fd{W_0}{u}\fd{W_0}{p}\Bigr)=0.
\end{equation*}

For $A^1$, one has
\begin{equation*}
  W_1=F^1p_0=\Bigl(p_3+\frac{2}{3}up_1+\frac{1}{3}u_1p_0\Bigr)p_0
  =p_3p_0+\frac{2}{3}up_1p_0.
\end{equation*}
Consequently,
\begin{equation*}
  \fd{W_1}{u}=\frac{2}{3}p_1p_0
\end{equation*}
and
\begin{multline*}
  \fd{W_1}{p}=\pd{W_1}{p_0}-D_x\pd{W_1}{p_1}-D_x^3\pd{W_1}{p_3}
  =-\Bigl(p_3+\frac{2}{3}up_1\Bigr)-D_x\Bigl(\frac{2}{3}up_0\Bigr)-D_x^3(p_0)\\
  =-2p_3-\frac{4}{3}up_1-\frac{2}{3}u_1p_0.
\end{multline*}
Hence,
\begin{equation*}
  \fd{W_1}{u}\fd{W_1}{p}
  =\frac{4}{3}p_0p_1p_3=D_x\Bigl(\frac{4}{3}p_0p_1p_2\Bigr)
\end{equation*}
that implies $\Eu(\fd{W_1}{u}\cdot\fd{W_1}{p})=0$.

\begin{remark}
  In~\cite{KerstenKrasilshchikVerbovetsky:InCSSREvDEq} we describe a
  class of equations which have the property that \eqref{eq:5}
  automatically implies the Hamiltonianity.  In particular, KdV
  belongs to that class, thus the above verification might be skipped.
\end{remark}

\subsubsection*{Nonlocal Hamiltonian structure}
Let us introduce a new (odd) nonlocal variable determined by the
equations
\begin{align*}
  r_x&=u_1p_0,\\
  r_t&=u_1p_2-u_2p_1+(uu_1+u_3)p_0
\end{align*}
(see Example~\ref{ex:cover}).  Then an additional solution of
equation~\eqref{eq:KdV-ham} arises:
\begin{equation*}
  F^2=p_5+\frac{4}{3}up_3+2u_1p_2+\Bigl(\frac{4}{9}u^2+\frac{4}{3}u_2\Bigr)p_1
  +\Bigl(\frac{4}{9}uu_1+\frac{1}{3}u_3\Bigr)p_0-\frac{1}{9}u_1r,
\end{equation*}
to which there corresponds the operator
\begin{equation*}
  A^2=D_x^5+\frac{4}{3}uD_x^3+2u_1D_x^2
  +\Bigl(\frac{4}{9}u^2+\frac{4}{3}u_2\Bigr)D_x
  +\Bigl(\frac{4}{9}uu_1+\frac{1}{3}u_3\Bigr)-\frac{1}{9}u_1D_x^{-1}\circ u_1.
\end{equation*}

\begin{remark}\label{rem:correspondence}
  Here and below we use the following correspondence between nonlocal
  variables and operators (in the case of evolution equations with
  one-dimensional~$x$). Let $p_k^j$ be the variables in the fibers of
  the $\ell_{\mathcal{E}}^*$-covering and nonlocal variable $r$ be
  determined by the relations
  \begin{equation*}
    \begin{cases}
      r_x=\sum_{k,j}a_k^jp_k^j,\\
      r_t=\sum_{k,j}b_k^jp_k^j
    \end{cases}
  \end{equation*}
  (cf.~Example~\ref{ex:cover}). Then the corresponding operator
  $\Delta_r$ acts on $\phi=(\phi^1,\dots,\phi^m)$ by
  \begin{equation*}
    \Delta_r(\phi)=D_x^{-1}\Big(\sum_{k,j}a_k^jD_x^k(\phi^j)\Big).
  \end{equation*}
\end{remark}

Simulating the techniques developed for the local case, it is a
straightforward check that $A^2$ is a Hamiltonian structure and all
three structures are pair-wise compatible. Moreover, they are related
to each other by the classical recursion operator
\begin{equation*}
  R=D_x^2+\frac{2}{3}u+\frac{1}{3}u_1D_x^{-1},
\end{equation*}
i.e., $A^1=R\circ A^0$ and $A^2=R\circ A^1$. In a similar way, one can
find a whole infinite series of nonlocal Hamiltonian structures for
the Korteweg-de Vries equation.

\begin{remark}\label{rem:simulate}
  We stress here the word \emph{simulating} above: at this moment, we
  do not have a consistent theory of Hamiltonian structures in the
  nonlocal setting. We hope to develop it elsewhere.
\end{remark}

\section{The Boussinesq equation}\label{sec:hamilt-struct-bouss}
In this section, we shall present, as another illustration of the
above developed methods, computation of local and nonlocal Hamiltonian
structures for the classical Boussinesq equation.  We consider this
equation as the system of the form
\begin{align}\label{eq:Bouss}
  u_t&=u_xv+uv_x+\sigma v_{xxx},\nonumber\\
  v_t&=u_x+vv_x,
\end{align}
where $\sigma\in\mathbb{R}$ is a constant.

All computations presented below were done by the software system
described in~\cite[Ch.~VIII]{KrasilshchikKersten:SROpCSDE} and we
expose here final results only.

\subsubsection*{The $\ell_{\mathcal{E}}^*$-covering}
The linearization operator restricted to $\mathcal{E}^\infty$ is
\begin{equation*}
  \ell_{\mathcal{E}}=\begin{pmatrix}
    vD_x-D_t&\sigma D_x^3+u_1\\
    D_x&vD_x+v_1-D_t
  \end{pmatrix}
\end{equation*}
while the adjoint one is expressed by
\begin{equation*}
  \ell_{\mathcal{E}}^*=\begin{pmatrix}
    -vD_x-v_1+D_t&-D_x\\
    -\sigma D_x^3+u_1&-vD_x+D_t
  \end{pmatrix}.
\end{equation*}
Hence, the $\ell^*_\mathcal{E}$-covering with the odd nonlocal
variables $p_i$, $q_i$ is defined by
\begin{align*}
  p_t&=vp_1+v_1p+q_1,\\
  q_t&=\sigma p_3-u_1p+vq_1.
\end{align*}

\subsubsection*{Local Hamiltonian operators}
In a completely similar way as described for the KdV equation in the
previous section, we solved the symmetry equation in the
$\ell^*_\mathcal{E}$-covering of the classical Boussinesq equation.
We found three local solutions of the form
\begin{align*}
  F^1 &= q_1,\\
  G^1 &= p_1;\\[1ex]
  F^2 &= 2\sigma p_3 + 2up_1 + u_1p_0 + vq_1,\\
  G^2 &= vp_1 + v_1p_0 + 2q_1 ;\\[1ex]
  F^3 &= 4\sigma vp_3 + 6\sigma v_1p_2+ 2(3\sigma v_2 + 2uv)p_1
  + 2(\sigma v_3 + u v_1 + u_1 v)p_0 \\
  &+ 4\sigma q_3 + (4u + v^2)q_1 + 2u_1q_0 ,\\
  G^3 &= 4\sigma p_3 + (4u + v^2)p_1 + 2(u_1 + v v_1)p_0 + 4vq_1 +
  2v_1q_0.
\end{align*}

In classical operator notation, they are represented as
\begin{equation*}
  A^1=
  \begin{pmatrix}
    0&D_x\\
    D_x&0
  \end{pmatrix},\quad
  A^2=
  \begin{pmatrix}
    2\sigma D_x^3 + 2uD_x + u_{1} &vD_x\\
    vD_x + v_{1}&2D_x
  \end{pmatrix},
\end{equation*}
while the third operator has the entries
\begin{align*}
  A_{11}^3 &= 4\sigma vD_x^3+ 6\sigma v_1 D_x^2
  + 2(3\sigma v_{2} + 2uv) D_x + 2 (\sigma v_3 + uv_1 + u_1v),\\
  A_{12}^3 &= 4 \sigma D_x^3 + (4 u + v^2) D_x + 2 u_{1} ,\\
  A_{21}^3 &= 4 \sigma D_x^3 + (4 u + v^2) D_x + 2 (u_{1} + v v_{1}),\\
  A_{22}^3 &= 4 v D_x + 2 v_{1}.
\end{align*}

\subsubsection*{Hamiltonianity and compatibility}
To test Hamiltonianity and compatibility conditions for the operators
$A^1$, $A^2$, $A^3$, we construct the bivectors
\begin{align*}
  W_{A^1}&=q_1p_0+p_1q_0,\\
  W_{A^2}&=(2\sigma p_3 + 2up_1 + u_1p_0 + vq_1)p_0+(vp_1 + v_1p_0 + 2q_1)q_0,\\
  W_{A^3}&=(4\sigma vp_3 + 6\sigma v_1p_2
  + 2(3\sigma v_2 + 2uv)p_1+ 2(\sigma v_3 + u v_1 + u_1 v)p_0 \\
  &+ 4\sigma q_3 + (4u + v^2)q_1 + 2u_1q_0)p_0\\
  &+(4\sigma p_3 + (4u + v^2)p_1 + 2(u_1 + v v_1)p_0 + 4vq_1 +
  2v_1q_0)q_0
\end{align*}
and straightforwardly check that
\begin{equation*}
  [\![W_{A^i},W_{A^j}]\!]=0,\qquad 1\le i\le j\le 3,
\end{equation*}
i.e., the operators $A^1$, $A^2$, $A^3$ meet both Hamiltonianity and
compatibility conditions.

\subsubsection*{Nonlocal Hamiltonian operators}
In order to describe nonlocal results we introduce three new nonlocal
variables $r_1$, $r_2$, $r_3$ over the $\ell^*_\mathcal{E}$-covering
by the following definitions
\begin{align*}
  r_{1,x} &= p_0u_1 + q_0v_1,\\
  r_{1,t} &= p_2\sigma v_1 - p_1\sigma v_2
  + p_0(\sigma v_3 + uv_1 + u_1v) + q_0(u_1 + vv_1);\\[1ex]
  r_{2,x} &= p_0(\sigma v_3 + uv_1 + u_1v) + q_0(u_1 + vv_1),\\
  r_{2,t} &= p_2\sigma (u_1 + vv_1) - p_1\sigma (u_2 + vv_2 + v_1^2) \\
  &+ p_0(\sigma u_3 + 2\sigma vv_3 + 3\sigma v_1v_2 + uu_1 + 2uvv_1 + u_1v^2) \\
  &+ q_0(\sigma v_3 + uv_1 + 2u_1 v + v^2v_1);\\[1ex]
  r_{3,x} &= p_0(4\sigma u_3 + 6 \sigma vv_3 + 12\sigma v_1v_2 + 6uu_1 + 6uvv_1 + 3u_1v^2) \\
  &+ q_0(4\sigma v_3 + 6uv_1 + 6u_1v + 3v^2 v_1),\\
  r_{3,t} &= p_2\sigma (4\sigma v_3 + 6uv_1 + 6u_1v + 3v^2v_1) \\
  &+ p_1\sigma (-4\sigma v_{4} - 6uv_2 - 12u_1v_1 - 6u_2v - 3v^2 v_2 - 6vv_1^2) \\
  &+ p_0(4\sigma^2 v_5 + 10\sigma uv_3 + 18\sigma u_1v_2
  + 18\sigma u_2v_1 + 10\sigma u_3v + 9 \sigma v^2 v_3 \\
  &+ 30 \sigma v v_1 v_2 + 6 \sigma v_1^3 + 6 u^2 v_1 + 12 u u_1 v + 9 u v^2 v_1 + 3 u_1 v^3) \\
  &+ q_0(4\sigma u_3 + 10\sigma vv_3 + 12\sigma v_1v_2 + 6uu_1 +
  12uvv_1 + 9u_1v^2 + 3v^3v_1).
\end{align*}
Using these nonlocal variables, we derived the following three
nonlocal Hamiltonian structures given by
\begin{align*}
  F^4 &= 8\sigma^2 p_5 + 2\sigma (8u + 3v^2)p_3
  + 6\sigma (4u_1 + 3vv_1) p_2\\
  &+ 2(8\sigma u_2 + 9\sigma vv_2 + 6\sigma v_1^2 + 4u^2 + 3uv^2) p_1\\
  &+ (4\sigma u_3 + 6\sigma vv_3 + 12\sigma v_1v_2 + 8uu_1 + 6uvv_1 +
  3u_1v^2)p_0
  + 12\sigma v q_3\\
  &+ 20\sigma v_1 q_2 + (16\sigma v_2 + 12uv + v^3) q_1
  + 2(2\sigma v_3 + 2uv_1 + 3u_1v) q_0\\
  &- 2u_1r_1,\\
  G^4 &= 12\sigma v p_3 + 16\sigma v_1 p_2
  + (12\sigma v_2 + 12uv + v^3) p_1\\
  &+ (4\sigma v_3 + 8uv_1 + 6u_1v + 3v^2v_1)p_0
  + 8\sigma q_3\\
  &+ 2(4u + 3v^2) q_1 + 2(2u_1 + 3vv_1) q_0
  - 2v_1r_1 ;\\[1ex]
  F^5 &= 32\sigma^2 v p_5 + 80\sigma^2 v_1 p_4
  + 8\sigma(14\sigma v_2 + 8uv + v^3) p_3\\
  &+ 4\sigma(22\sigma v_3 + 24uv_1 + 24u_1v + 9v^2v_1) p_2
  + 4(10 \sigma^2 v_4 + 20\sigma uv_2 \\
  &+ 26\sigma u_1v_1 + 16\sigma u_2v + 9\sigma v^2v_2 + 12\sigma vv_1^2 + 8u^2 v + 2 u v^3) p_1\\
  &+ 4(2\sigma^2 v_5 + 6\sigma uv_3 + 11\sigma u_1v_2
  + 9\sigma u_2v_1 + 4\sigma u_3v + 3\sigma v^2v_3 \\
  &+ 12\sigma vv_1 v_2 + 3\sigma v_1^3 + 4 u^2 v_1 + 8uu_1v + 3uv^2v_1
  + u_1v^3)p_0
  + 16\sigma^2 q_5\\
  &+ 8\sigma (4 u + 3 v^2) q_3 + 16\sigma(3u_1 + 5vv_1) q_2
  + (32\sigma u_2 + 64\sigma vv_2 + 44\sigma v_1^2 \\
  &+ 16u^2 + 24uv^2 + v^4) q_1
  + 4(2\sigma u_3 + 4\sigma vv_3 + 6\sigma v_1v_2 + 4 uu_1 \\
  &+ 4uvv_1 + 3u_1v^2) q_0 - 4 u_1 r_2
  - 4(\sigma v_{3} + uv_1 + u_1v) r_1,\\
  G^5 &= 16\sigma^2 p_5 + 8\sigma (4u + 3v^2) p_3
  + 16\sigma (3u_1 + 4v v_1) p_2\\
  &+ (32\sigma u_2 + 48\sigma vv_2 + 28\sigma v_1^2 + 16 u^2 + 24uv^2 + v^4) p_1\\
  &+ 4(2\sigma u_3 + 4\sigma vv_3 + 8\sigma v_1v_2 + 4uu_1 + 8uvv_1 + 3u_1v^2 + v^3v_1) p_0\\
  &+ 32\sigma v q_3 + 48\sigma v_1 q_2
  + 8(4\sigma v_2 + 4uv + v^3) q_1\\
  &+ 4(2\sigma v_3 + 4uv_1 + 4u_1v + 3v^2v_1) q_0 - 4v_1 r_2
  - 4(u_1 + vv_1) r_1;\\[1ex]
  F^6 &= - 32\sigma^3 p_7 - 16\sigma^2 (6u + 5 v^2) p_5
  - 80\sigma^2(3u_1 + 5 vv_1) p_4\\
  &- 2\sigma(160\sigma u_2 + 280\sigma vv_2 + 204\sigma v_1^2 + 48u^2 + 80uv^2 + 5v^4) p_3\\
  &- 4\sigma (60\sigma u_3 + 110\sigma vv_3 + 216\sigma v_1v_2 + 72uu_1 + 120uvv_1 + 60u_1v^2 \\
  &+ 15v^3v_1) p_2
  - 2(48\sigma^2 u_4 + 100\sigma^2 vv_4 + 244\sigma^2 v_1v_3 + 168\sigma^2 v_2^2 \\
  &+ 96\sigma uu_2 + 200\sigma uvv_2 + 136\sigma uv_1^2
  + 68\sigma u_1^2 + 260\sigma u_1vv_1 + 80\sigma u_2v^2 \\
  &+ 30\sigma v^3v_2 + 60\sigma v^2v_1^2 + 16u^3 + 40u^2v^2 + 5 uv^4) p_1\\
  &- (16\sigma^2 u_5 + 40\sigma^2 vv_5 + 120\sigma^2 v_1v_4 + 208\sigma^2 v_2v_3 + 48\sigma uu_3 \\
  &+ 120\sigma uvv_3 + 232\sigma uv_1v_2 + 88\sigma u_1u_2 + 220\sigma u_1vv_2 + 156\sigma u_1v_1^2 \\
  &+ 180\sigma u_2vv_1 + 40\sigma u_3v^2 + 20\sigma v^3v_3 + 120\sigma v^2v_1v_2 + 60\sigma vv_1^3 \\
  &+ 48u^2u_1 + 80u^2vv_1 + 80uu_1v^2 + 20uv^3v_1 + 5u_1v^4) p_0
  - 80\sigma^2 v q_5\\
  &- 224\sigma^2 v_1 q_4 - 40\sigma (8\sigma v_2 + 4uv + v^3) q_3
  - 8\sigma (30\sigma v_3 + 32uv_1 \\
  &+ 30u_1v + 25v^2 v_1) q_2
  - (96\sigma^2 v_4 + 192\sigma uv_2 + 256\sigma u_1v_1 \\
  &+ 160\sigma u_2v + 160\sigma v^2v_2 + 220\sigma vv_1^2 + 80 u^2 v + 40 uv^3 + v^5) q_1\\
  &- 4(4\sigma^2 v_5 + 12\sigma uv_3 + 22\sigma u_1 v_2 + 18\sigma u_2v_1 + 10\sigma u_3v + 10\sigma v^2v_3 \\
  &+ 30\sigma vv_1v_2 + 6\sigma v_1^3 + 8u^2v_1 + 20uu_1v + 10uv^2v_1 + 5u_1v^3) q_0\\
  &+ 2u_1 r_3
  + 8(\sigma v_3 + uv_1 + u_1v) r_2\\
  &+ 2(4\sigma u_3 + 6\sigma vv_3 + 12\sigma v_1v_2 + 6uu_1 + 6uvv_1 + 3u_1v^2)r_1,\\
  G^6 &= - 80\sigma^2 v p_5 - 176\sigma^2 v_1 p_4
  - 8\sigma(28\sigma v_2 + 20uv - 5v^3) p_3\\
  &- 16\sigma (11\sigma v_3 + 14uv_1 + 15u_1v + 10v^2v_1) p_2\\
  &- (80\sigma^2 v_4 + 160\sigma uv_2 + 224\sigma u_1v_1 + 160\sigma u_2v + 120\sigma v^2v_2 \\
  &+ 140\sigma vv_1^2 + 80u^2v + 40uv^3 + v^5) p_1
  - (16\sigma^2 v_5 + 48\sigma uv_3 \\
  &+ 88\sigma u_1v_2 + 88\sigma u_2v_1 + 40\sigma u_3v + 40\sigma v^2v_3 + 160\sigma vv_1v_2 + 36\sigma v_1^3 \\
  &+ 48u^2v_1 + 80uu_1v + 80uv^2v_1 + 20u_1v^3 + 5v^4v_1) p_0
  - 32\sigma^2 q_5\\
  &- 16\sigma (4u + 5v^2) q_3
  - 48\sigma (2u_1 + 5 v v_1) q_2\\
  &- 2(32\sigma u_2 + 80\sigma vv_2 + 52\sigma v_1^2 + 16u^2 + 40uv^2 + 5v^4) q_1\\
  &- 4(4\sigma u_3 + 10\sigma vv_3 + 16\sigma v_1v_2 + 8uu_1 + 20uvv_1 + 10u_1v^2 \\
  &+ 5v^3v_1) q_0 + 2v_1 r_3
  + 8(u_1 + vv_1) r_2\\
  &+ 2(4\sigma v_3 + 6 uv_1 + 6u_1v + 3v^2v_1) r_1.
\end{align*}

In the conventional notation the operators $A^4$, $A^5$, and $A^6$
have the following entries:
\begin{align*}
  A^4_{11} &= 8 \sigma^2 D_x^5 + 2 \sigma (8 u + 3 v^2) D_x^3
  + 6 \sigma (4 u_{1} + 3 v v_{1}) D_x^2 \\
  &+ 2 (8 \sigma u_{2} + 9 \sigma v v_{2} + 6 \sigma v_{1}^2 + 4 u^2 + 3 u v^2) D_x^1 \\
  &+ (4 \sigma u_{3} + 6 \sigma v v_{3} + 12 \sigma v_{1} v_{2} + 8 u u_{1} + 6 u v v_{1} + 3 u_{1} v^2)\\
  &- 2 u_{1} D_x^{-1}\circ u_{1},\\
  A^4_{12} &= 12 \sigma v D_x^3 + 20 \sigma v_{1} D_x^2
  + (16 \sigma v_{2} +12 u v + v^3) D_x \\
  &+ 2 (2 \sigma v_{3} + 2 u v_{1} + 3 u_{1} v)
  - 2 u_{1} D_x^{-1}\circ v_{1},\\
  A^4_{21} &=12 \sigma v D_x^3 + 16 \sigma v_{1} D_x^2
  + (12 \sigma v_{2} + 12 u v + v^3) D_x \\
  &+ (4 \sigma v_{3} + 8 u v_{1} + 6 u_{1} v + 3 v^2 v_{1})
  - 2 v_{1} D_x^{-1}\circ u_{1},\\
  A^4_{22} &= 8 \sigma D_x^3 + 2 (4 u + 3 v^2)D_x + 2 (2 u_{1} + 3 v
  v_{1})- 2 v_{1} D_x^{-1}\circ v_{1}.
\end{align*}
The matrix elements of $A^5$ are given as
\begin{align*}
  A^5_{11} &= 32 \sigma^2 v D_x^5 + 80 \sigma^2 v_{1} D_x^4
  + 8 \sigma (14 \sigma v_{2} + 8 u v + v^3) D_x^3 \\
  &+ 4 \sigma (22 \sigma v_{3} + 24 u v_{1} + 24 u_{1} v + 9v^2 v_{1}) D_x^2 \\
  &+ 4 (10 \sigma^2 v_{4} + 20 \sigma u v_{2} + 26 \sigma u_{1}v_{1} + 16 \sigma u_{2} v + 9 \sigma v^2 v_{2} + 12 \sigma v v_{1}^2 \\
  &+ 8u^2 v + 2 u v^3) D_x
  + 4 (2 \sigma^2 v_{5} + 6 \sigma u v_{3} + 11 \sigma u_{1}v_{2} + 9 \sigma u_{2} v_{1} \\
  &+ 4 \sigma u_{3} v + 3 \sigma v^2 v_{3} + 12 \sigma v v_{1} v_{2} + 3 \sigma v_{1}^3 + 4 u^2 v_{1} + 8 u u_{1} v \\
  &+ 3 u v^2 v_{1} + u_{1} v^3)\\
  &- 4 (\sigma v_{3} + u v_{1} + u_{1} v)D_x^{-1}\circ u_{1}
  - 4 u_{1} D_x^{-1}\circ (\sigma v_{3} + uv_{1} + u_{1}v),\\
  A^5_{12} &=16 \sigma^2 D_x^5 + 8 \sigma (4 u + 3 v^2) D_x^3
  + 16 \sigma (3 u_{1} + 5v v_{1}) D_x^2 \\
  &+ (32 \sigma u_{2} + 64 \sigma v v_{2} + 44 \sigma v_{1}^2 + 16 u^2 + 24 u v^2 + v^4) D_x \\
  &+ 4 (2 \sigma u_{3} + 4 \sigma v v_{3} + 6 \sigma v_{1} v_{2} +
  4 u u_{1} + 4 u v v_{1} + 3 u_{1} v^2)\\
  &- 4(\sigma v_{3} + u v_{1} + u_{1} v) D_x^{-1}\circ v_{1}
  - 4 u_{1} D_x^{-1}\circ (u_{1} + vv_{1}), \\
  A^5_{21} &= + 16 \sigma^2 D_x^5+ 8 \sigma (4 u + 3 v^2) D_x^3
  + 16 \sigma (3 u_{1} + 4 v v_{1}) D_x^2\\
  &+ (32 \sigma u_{2} + 48 \sigma v v_{2} + 28 \sigma v_{1}^2 + 16 u^2 + 24 u v^2 + v^4) D_x \\
  &+ 4 (2 \sigma u_{3} + 4 \sigma v v_{3} + 8 \sigma v_{1} v_{2} + 4 u u_{1} + 8 u v v_{1} + 3 u_{1} v^2 + v^3 v_{1})\\
  &- 4 (u_{1} + v v_{1})D_x^{-1}\circ u_{1}
  - 4 v_{1} D_x^{-1}\circ (\sigma v_{3} + uv_{1} + u_{1}v) ,\\
  A^5_{22} &= 32 \sigma v D_x^3 + 48 \sigma v_{1} D_x^2
  + 8 (4 \sigma v_{2} + 4 u v + v^3) D_x \\
  &+ 4 (2 \sigma v_{3} + 4 u v_{1} + 4 u_{1} v + 3 v^2 v_{1})\\
  &- 4 (u_{1} + v v_{1}) D_x^{-1}\circ v_{1} - 4 v_{1} D_x^{-1}\circ
  (u_{1} + vv_{1}).
\end{align*}
The matrix elements of $A^6$ are given as
\begin{align*}
  A^6_{11} &= - 32 \sigma^3 D_x^7 + 16 \sigma^2 ( - 6 u - 5 v^2) D_x^5
  + 80 \sigma^2 ( - 3 u_{1} - 5 v v_{1}) D_x^4 \\
  &+ 2 \sigma ( - 160 \sigma u_{2} - 280 \sigma v v_{2} - 204 \sigma v_{1}^2 - 48 u^2 - 80 u v^2 - 5 v^4) D_x^3\\
  &+ 4 \sigma ( - 60 \sigma u_{3} - 110 \sigma v v_{3} - 216\sigma v_{1} v_{2} - 72 u u_{1} - 120 u v v_{1} \\
  &- 60 u_{1} v^2 - 15 v^3v_{1}) D_x^2
  + 2 ( - 48 \sigma^2 u_{4} - 100 \sigma^2 v v_{4} - 244\sigma^2 v_{1} v_{3} \\
  &- 168 \sigma^2 v_{2}^2 - 96 \sigma u u_{2} - 200\sigma u v v_{2}- 136 \sigma u v_{1}^2 - 68 \sigma u_{1}^2 \\
  &- 260\sigma u_{1} v v_{1} - 80 \sigma u_{2} v^2 - 30 \sigma v^3 v_{2} - 60\sigma v^2 v_{1}^2 - 16 u^3 \\
  &- 40 u^2 v^2 - 5 u v^4) D_x \\
  &+ ( - 16 \sigma^2 u_{5} - 40 \sigma^2 v v_{5} - 120\sigma^2 v_{1} v_{4} - 208 \sigma^2 v_{2} v_{3} - 48 \sigma u u_{3} \\
  &- 120 \sigma u v v_{3} - 232 \sigma u v_{1} v_{2} - 88 \sigma u_{1} u_{2} - 220 \sigma u_{1} v v_{2} - 156 \sigma u_{1} v_{1}^2 \\
  &- 180 \sigma u_{2} v v_{1} - 40 \sigma u_{3} v^2 - 20 \sigma v^3 v_{3} - 120 \sigma v^2 v_{1} v_{2} - 60 \sigma v v_{1}^3 \\
  &- 48 u^2 u_{1} - 80 u^2 v v_{1} - 80 u u_{1} v^2 - 20 u v^3 v_{1} - 5 u_{1} v^4)\\
  &+ 2(4 \sigma u_{3} + 6 \sigma v v_{3} + 12 \sigma v_{1} v_{2} +6 u u_{1} + 6 u v v_{1} + 3 u_{1} v^2) D_x^{-1}\circ u_{1} \\
  &+ 8 (\sigma v_{3} + u v_{1} + u_{1} v) D_x^{-1}\circ (\sigma v_{3} + uv_{1} + u_{1}v) \\
  &+ 2 u_{1} D_x^{-1}\circ (4 \sigma u_{3} + 6 \sigma v v_{3} + 12 \sigma v_{1} v_{2} + 6 u u_{1} + 6 u v v_{1} + 3 u_{1} v^2) ,\\
  A^6_{12} &= - 80 \sigma^2 v D_x^5 - 224 \sigma^2 v_{1} D_x^4
  + 40 \sigma ( - 8 \sigma v_{2} - 4 u v - v^3) D_x^3 \\
  &+ 8 \sigma ( - 30 \sigma v_{3} - 32 u v_{1} - 30 u_{1} v - 25 v^2 v_{1}) D_x^2 \\
  &+ ( - 96 \sigma^2 v_{4} - 192 \sigma u v_{2} - 256 \sigma u_{1} v_{1} - 160 \sigma u_{2} v - 160 \sigma v^2 v_{2} \\
  &- 220 \sigma vv_{1}^2 - 80 u^2 v - 40 u v^3 - v^5) D_x
  + 4 ( - 4 \sigma^2 v_{5} - 12 \sigma u v_{3} \\
  &- 22 \sigma u_{1} v_{2} - 18 \sigma u_{2} v_{1} - 10 \sigma u_{3} v - 10 \sigma v^2 v_{3}- 30 \sigma v v_{1} v_{2} - 6 \sigma v_{1}^3 \\
  &- 8 u^2 v_{1} - 20 u u_{1} v - 10u v^2 v_{1} - 5 u_{1} v^3)\\
  &+2 (4 \sigma u_{3} + 6 \sigma v v_{3} + 12 \sigma v_{1} v_{2} + 6 u u_{1} + 6 u v v_{1} + 3 u_{1} v^2) D_x^{-1}\circ v_{1} \\
  &+ 8 (\sigma v_{3} + u v_{1} + u_{1} v)D_x^{-1}\circ (u_{1} + vv_{1}) \\
  &+ 2 u_{1} D_x^{-1}\circ (4 \sigma v_{3} + 6 u v_{1} + 6 u_{1} v + 3 v^2 v_{1}) ,\\
  A^6_{21} &= - 80 \sigma^2 v D_x^5 - 176 \sigma^2 v_{1} D_x^4
  + 8 \sigma ( - 28 \sigma v_{2} - 20 u v - 5 v^3) D_x^3 \\
  &+ 16\sigma ( - 11 \sigma v_{3} - 14 u v_{1} - 15 u_{1} v - 10 v^2 v_{1}) D_x^2 \\
  &+ ( - 80 \sigma^2 v_{4} - 160 \sigma u v_{2} - 224 \sigma u_{1} v_{1} - 160 \sigma u_{2} v - 120 \sigma v^2 v_{2} \\
  &- 140 \sigma v v_{1}^2 - 80 u^2 v - 40 u v^3 - v^5) D_x
  + ( - 16 \sigma^2 v_{5} - 48 \sigma u v_{3} \\
  &- 88 \sigma u_{1} v_{2} - 88 \sigma u_{2} v_{1} - 40 \sigma u_{3} v - 40 \sigma v^2 v_{3} - 160 \sigma v v_{1} v_{2} - 36 \sigma v_{1}^3 \\
  &- 48u^2 v_{1} - 80 u u_{1} v - 80 u v^2 v_{1} - 20 u_{1} v^3 - 5 v^4 v_{1})\\
  &+ 2 (4 \sigma v_{3} + 6 u v_{1} + 6 u_{1} v + 3 v^2 v_{1})D_x^{-1}\circ u_{1} \\
  &+ 8 (u_{1} + v v_{1})D_x^{-1}\circ (\sigma v_{3} + uv_{1} + u_{1}v) \\
  &+ 2 v_{1} D_x^{-1}\circ (4 \sigma u_{3} + 6 \sigma v v_{3} + 12 \sigma v_{1} v_{2} + 6 u u_{1} + 6 u v v_{1} + 3 u_{1} v^2) ,\\
  A^6_{22} &=- 32 \sigma^2 D_x^5 + 16 \sigma ( - 4 u - 5 v^2) D_x^3
  + 48 \sigma ( - 2 u_{1} - 5 v v_{1}) D_x^2 \\
  &+ 2 ( - 32 \sigma u_{2} - 80 \sigma v v_{2} - 52 \sigma v_{1}^2 - 16 u^2 - 40 u v^2 - 5 v^4) D_x \\
  &+ 4 ( - 4 \sigma u_{3} - 10 \sigma v v_{3} - 16 \sigma v_{1}v_{2} - 8 u u_{1} - 20 u v v_{1} \\
  &- 10 u_{1} v^2 - 5 v^3 v_{1})\\
  &+ 2(4 \sigma v_{3} + 6 u v_{1} + 6 u_{1} v + 3 v^2 v_{1}) D_x^{-1}\circ v_{1} \\
  &+ 8 (u_{1} + v v_{1})D_x^{-1}\circ (u_{1} + vv_{1}) \\
  &+ 2 v_{1} D_x^{-1}\circ (4 \sigma v_{3} + 6 u v_{1} + 6 u_{1} v + 3
  v^2 v_{1}).
\end{align*}

Similar to the previous cases, we checked the conditions for
Hamiltonianity and compatibility of all six Hamiltonian structures. It
is also easy to check that all six structures are related by the
recursion operator constructed for symmetries of the Boussinesq
equation, see, e.g.,~\cite{KrasilshchikKersten:SROpCSDE}.

\section{The coupled KdV-mKdV system}\label{sec:appl}

We shall now describe a Hamiltonian structure for the coupled KdV-mKdV
system of the form
\begin{align}\label{eq:kdv-mkdv}
  u_t &= - u_{xxx} + 6uu_x - 3vv_{xxx} - 3v_xv_{xx} + 3u_xv^2 + 6uvv_x,\nonumber\\
  v_t &= - v_{xxx} + 3v^2v_x + 3uv_x + 3u_xv.
\end{align}
This system arises as the so-called \emph{bosonic limit} of the $N=2$,
$a=1$ supersymmetric extension of the KdV
equation~\cite{Kersten:SROpSKEq}; integrability properties of this
system (existence of a recursion operator) were studied
in~\cite{KerstenKrasilshchik:CInCKS}.
In~\cite{Karasu(Kalkanli)SakovichYurdusen:InKKCKEqSAnLP}, by means of
the prolongations structure techniques, a Lax pair
for~\eqref{eq:kdv-mkdv} was constructed.

Denote the evolution equation corresponding to~\eqref{eq:kdv-mkdv} by
$\mathcal{E}^\infty$ and choose for coordinates in
$\mathcal{E}^\infty$ the functions
\begin{equation*}
  x,\ t,\ u=u_0,\ v=v_0,\dots,u_k,\ v_k,\dots
\end{equation*}
Then the total derivative operators restricted to $\mathcal{E}^\infty$
are written in the form
\begin{align}\label{eq:KMK-tot}
  D_x&=\pd{}{x}+\sum_{k\ge0}\Bigl(u_{k+1}\pd{}{u_k}
  +v_{k+1}\pd{}{v_k}\Bigr),\nonumber\\
  D_t&=\pd{}{t}+\sum_{k\ge0}\Bigl(D_x^k(f)\pd{}{u_k}
  +D_x^k(g)\pd{}{v_k}\Bigr),
\end{align}
where
\begin{align*}
  f &= - u_3 + 6uu_1 - 3vv_3 - 3v_1v_2 + 3u_1v^2 + 6uvv_1,\\
  g &= - v_3 + 3v^2v_1 + 3uv_1 + 3u_1v
\end{align*}
are the functions at the right-hand side of~\eqref{eq:kdv-mkdv}.

\subsubsection*{The $\ell_{\mathcal{E}}^*$-covering}
The linearization operator restricted to $\mathcal{E}^\infty$ is
\begin{equation}\label{eq::KMK-ell}
  \ell_{\mathcal{E}}=
  \begin{pmatrix}
    \ell_{11}&\ell_{12}\\
    \ell_{21}&\ell_{22}
  \end{pmatrix},
\end{equation}
where
\begin{align*}
  \ell_{11}&=D_t+D_x^3-(6u+3v^2)D_x-6(u_1+vv_1),\\
  \ell_{12}&=3vD_x^3+3v_1D_x^2-(6uv-3v_2)D_x-6u_1v-6uv_1+3v_3,\\
  \ell_{21}&=-3vD_x-3v_1,\\
  \ell_{22}&=D_t+D_x^3-3(u+v^2)D_x-(6vv_1+3u).
\end{align*}
Consequently, the adjoint operator is
\begin{equation}\label{eq::KMK-ell*}
  \ell_{\mathcal{E}}^*=
  \begin{pmatrix}
    \ell_{11}^*&\ell_{21}^*\\
    \ell_{12}^*&\ell_{22}^*
  \end{pmatrix},
\end{equation}
where
\begin{align*}
  \ell_{11}^*&=-D_t-D_x^3+(6u+3v^2)D_x,\\
  \ell_{21}^*&=3vD_x,\\
  \ell_{12}^*&=-3vD_x^3-6v_1D_x^2+6(uv-v_2)D_x,\\
  \ell_{22}^*&=-D_t-D_x^3+3(u+v^2)D_x.
\end{align*}

Following the general theory of Section~\ref{sec:hamilt-evol-equat},
we now construct the $\ell^*_\mathcal{E}$-covering for the equation
$\mathcal{E}^\infty$ by introducing new odd variables $p=p_0$,
$q=q_0,\dots,p_k$, $q_k,\dots$, $p_k=D_x^k(p)$, $q_k=D_x^k(q)$, that
obey the equations
\begin{align}
  p_t&=-p_3+(6u+3v^2)p_1 + 3vq_1,\label{ell*1}\\
  q_t&=-3vp_3-6v_1p_2+6(uv-v_2)p_1-q_3+3(u+v^2)q_1.\label{ell*2}
\end{align}

\subsubsection*{Solving the defining equations} We now introduce a vector-function of the form
\begin{equation*}
  \begin{pmatrix}
    F\\
    G
  \end{pmatrix}
  =
  \begin{pmatrix}
    \sum_i(F_i^up_i+F_i^vq_i)\\
    \sum_i(G_i^up_i+G_i^vq_i)
  \end{pmatrix},
\end{equation*}
where $F_i^u$, $F_i^v$, $G_i^u$, $G_i^v$ are functions on
$\mathcal{E}^\infty$, and solve the equation
\begin{equation}\label{eq:ell*ell}
  \begin{pmatrix}
    \ell_{11}&\ell_{12}\\
    \ell_{21}&\ell_{22}
  \end{pmatrix}
  \begin{pmatrix}
    F\\
    G
  \end{pmatrix}
  =0.
\end{equation}
The operators $\ell_{ij}$ here are lifted to the
$\ell^*_\mathcal{E}$-covering, which means that the total derivatives
are now of the form
\begin{align}\label{eq:ell*-tot}
  \tilde{D}_x&=\pd{}{x}+\sum_{k\ge0}\Bigl(u_{k+1}\pd{}{u_k}
  +v_{k+1}\pd{}{v_k}+p_{k+1}\pd{}{p_k}+q_{k+1}\pd{}{q_k}\Bigr),
  \nonumber\\
  \tilde{D}_t&=\pd{}{t}+\sum_{k\ge0}\Bigl(D_x^k(f)\pd{}{u_k}
  +D_x^k(g)\pd{}{v_k}
  +D_x^k(f')\pd{}{p_k}+D_x^k(g')\pd{}{q_k}\Bigr),
\end{align}
where $f'$ and $g'$ are the right-hand sides of~\eqref{ell*1}
and~\eqref{ell*2}, respectively.

The following solution was obtained
\begin{align*}
  F&= -p_3 + 4up_1 + 2u_1p_0 + 2vq_1,\\
  G&= 2vp_1 + 2v_1p_0 + q_1,
\end{align*}
to which there corresponds the operator
\begin{equation}\label{eq:Ham3}
  A=
  \begin{pmatrix}
    -D_x^3 + 4uD_x + 2u_1&2vD_x\\
    2vD_x + 2v_1 &D_x
  \end{pmatrix}.
\end{equation}

\subsubsection*{The Hamiltonianity test}
We shall check now that the operator $A$ presented by~\eqref{eq:Ham3}
is Hamiltonian.  The first property is obvious: evidently, $A^*=-A$,
i.e., $A$ is a skew-adjoint operator.

To check the second property, we construct the bivector
\begin{multline*}
  W_A=Fp_0+Gq_0=(-p_3 + 4up_1 + 2u_1p_0 + 2vq_1)p_0+(2vp_1 + 2v_1p_0 + q_1)q_0\\
  =p_0p_3 - 4up_0p_1 - 2vp_0q_1 + 2vp_1q_0 + 2v_1p_0q_0 - q_0q_1
\end{multline*}
and verify condition~\eqref{eq:2}, i.e.,
\begin{equation}\label{eq:Ham-cond}
  \Eu\Bigl(\fd{W_A}{u}\fd{W_A}{p}+\fd{W_A}{v}\fd{W_A}{q}\Bigr)=0.
\end{equation}
But
\begin{align*}
  \fd{W_A}{u}&= -4p_0p_1,\\
  \fd{W_A}{v}&= -4p_0 q_1,\\
  \fd{W_A}{p}&= 2(p_3 - 4up_1 - 2u_1p_0 - 2vq_1),\\
  \fd{W_A}{q}&= 2(- 2vp_1 - 2v_1p_0 - q_1)
\end{align*}
and consequently
\begin{equation*}
  \fd{W_A}{u}\fd{W_A}{p}+\fd{W_A}{v}\fd{W_A}{q} = - 8 p_0 p_1 p_3 = D_x(- 8 p_0 p_1 p_2),
\end{equation*}
i.e.,~\eqref{eq:Ham-cond} holds.

\subsubsection*{Existence of a Hamiltonian}
Let us show that the KdV-mKdV system~\eqref{eq:kdv-mkdv} possesses a
Hamiltonian, i.e., its right-hand side may be represented in the form
\begin{equation}\label{eq:Ham-repr}
  \begin{pmatrix}
    f\\
    g
  \end{pmatrix}
  =A\Eu(X),
\end{equation}
where $A$ is the Hamiltonian operator described above and $X$ is the
$dx$-component of a conservation law $\eta=X\,dx+T\,dt$ (the energy).

We computed directly several conservation laws of lower order and
obtained the following results (for the sake of briefness, we omit the
corresponding $dt$-components):
\begin{align*}
  \eta_1:\quad X&=  v,\\
  \eta_2:\quad X&=  u,\\
  \eta_4:\quad X&=  \frac{1}{2}(u^2 + u v^2 - v v_2), \\
  \eta_6:\quad X&= 12 u^3 + 24u^2v^2-6uu_2+6uv^4-30uvv_2-3u_2 v^2 -8
  v^3 v_2 + 6 v v_4.
\end{align*}

Generating functions corresponding to these conservation laws, that
is, vector-functions of the form
\begin{equation*}
  \begin{pmatrix}
    \phi\\
    \psi
  \end{pmatrix}
  =\Eu(X)=
  \begin{pmatrix}
    \fd{X}{u}\\
    \fd{X}{v}
  \end{pmatrix},
\end{equation*}
are
\begin{align*}
  \phi_1&=0,&\psi_1&=1;\\
  \phi_2&=1,&\psi_2&=0;\\
  \phi_4&=u + \frac{1}{2}v^2,&\psi_4&=uv - v_2;\\
  \phi_6&=6(6u^2 + 8uv^2 - 2u_2 + v^4
  &\psi_6&=12(4u^2v + 2uv^3 - 5uv_2 - 5u_1v_1 \\
  &- 6vv_2 - v_1^2);&&- 3u_2v- 4v^2v_2 - 4vv_1^2 + v_4).
\end{align*}
Applying $A$ to $\Eu(X)$, where $X$ corresponds to $\eta_4$, we see
that~\eqref{eq:Ham-repr} holds.
\begin{theorem}\label{thm:Ham-kdv-mkdv}
  The coupled KdV-mKdV system~\eqref{eq:kdv-mkdv} is Hamiltonian with
  respect to the the Hamiltonian operator~\eqref{eq:Ham3} and
  possesses the Hamiltonian $X=\frac{1}{2}(u^2 + u v^2 - v v_2)$. The
  corresponding energy is given by the form
  \begin{multline*}
    \eta=\frac{1}{2}(u^2 + u v^2 - v v_2)\,dx 
    +\frac{1}{2}(4 u^3+9u^2v^2-2uu_2+3uv^4-11uvv_2+uv_1^2\\
    +u_1^2-u_1vv_1-4u_2v^2-6v^3v_2-3v^2v_1^2+vv_4-v_1v_3+v_2^2)\,dt.
  \end{multline*}
  This structure is unique in the class of Hamiltonian structures
  independent of $x$ and $t$ and polynomial in $u_k$\textup{,} $v_k$.
\end{theorem}

\begin{proof}[Proof of uniqueness]
  Let us first note that equation~\eqref{eq:kdv-mkdv} admits a scaling
  symmetry that allows to assign \emph{gradings} to all variables $x$,
  $t$, $u_k$, and $v_k$:
  \begin{equation*}
    \abs{x}=-1,\ \abs{t}=-3,\ \abs{u_k}=k+2,\ \abs{v_k}=k+1;
  \end{equation*}
  the grading of a monomial is the sum of gradings of the factors
  entering this monomial. In particular, $\abs{f}=5$, $\abs{g}=4$. All
  constructions are in agreement with these gradings and we may
  restrict computations to homogeneous components. Since the grading
  of the expression $B\Eu(X)$ is $\abs{B}+\abs{\Eu(X)}$, we conclude
  that the grading of the generating function $\Eu(X)$ is less than
  that of the right-hand side of~\eqref{eq:kdv-mkdv}. This fact
  restricts the choice of possible Hamiltonians just to several ones
  and by a direct computation we find that the only possible solution
  is given by Theorem~\ref{thm:Ham-kdv-mkdv}.
\end{proof}

\subsubsection*{Discussion\textup{:} nonlocalities}
In spite of the previous result, we have constructed another
Hamiltonian operator for the system under consideration. This operator
exists in an appropriate nonlocal setting. First, we introduce a new
nonlocal variable $w$ defined by
\begin{align*}
  w_x&=v,\\
  w_t&=3uv+v^3-v_2
\end{align*}
and corresponding to the conservation law $\eta_1$ (see
Example~\ref{ex:cover}).

In this nonlocal setting, it is possible to extend the
$\ell^*_\mathcal{E}$-covering by adding odd nonlocal variables $r_1$,
$r_2$, $r_3$ defined by the relations
\begin{align*}
  r_{1,x}&= q_0v_1 + p_0 u_1,\\
  r_{1,t} &= - q_2v_1 + p_2(- u_1 - 3vv_1) + q_1v_2 + p_1(u_2 + 3vv_2 - 3v_1^2) \\
  &+ q_0(3uv_1 + 3u_1v + 3v^2v_1 - v_3) \\
  &+ p_0(6uu_1 + 6uvv_1 + 3u_1v^2 - u_3 - 3vv_3 - 3v_1v_2);\\[1ex]
  r_{2,x} &= \frac{1}{2}q_0 \cos(2w)v_1 - \frac{1}{2}q_0\sin(2w)u \\
  &- p_0\cos(w)(\frac{1}{2}u_1+vv_1)+ p_0 \sin(2w)(uv - \frac{1}{2}v_2),\\
  r_{2,t} &= -\frac{1}{2}q_2\cos(2w)v_1 + \frac{1}{2}q_2\sin(2w)u + \frac{1}{2}p_2\cos(2w)(u_1 - vv_1)\\
  &+ \frac{1}{2}p_2\sin(2w)(uv + v_2) - q_1\cos(2w)(uv + \frac{1}{2}v_2) \\
  &- q_1\sin(2w)(\frac{1}{2}u_1 + vv_1) \\
  &- p_1\cos(2w)(uv^2 + \frac{1}{2}u_2 + \frac{1}{2}vv_2 + \frac{5}{2}v_1^2) \\
  &+ p_1\sin(2w)(\frac{5}{2}uv_1 + \frac{1}{2}u_1v - v^2v_1 - \frac{1}{2}v_3) \\
  &+ \frac{1}{2}q_0\cos(2w)(5uv_1 + u_1v + v^2v_1 - v_3) \\
  &+ q_0\sin(2w)(-\frac{3}{2}u^2 - \frac{1}{2}uv^2 + \frac{1}{2}u2 + \frac{1}{2}vv_2 + v_1^2) \\
  &+ p_0\cos(2w)(-3uu_1 - 6uvv_1 - \frac{3}{2}u_1v^2 + \frac{1}{2}u_3 - v^3v_1 \\
  &+ \frac{3}{2}vv_3 + \frac{5}{2}v_1v_2) \\
  &+ p_0\sin(2w)(3u^2v + uv^3 - \frac{5}{2}uv_2 - 3u_1v_1 - \frac{3}{2}u_2v - \frac{3}{2}v^2v_2 \\
  &- 3vv_1^2 + \frac{1}{2}v_4);\\[1ex]
  r_{3,x} &= -\frac{1}{2}q_0\cos(2w)u - \frac{1}{2}q_0\sin(2w)v_1 \\
  &+ p_0\cos(2w)(uv - \frac{1}{2}v_2) + p_0\sin(2w)(\frac{1}{2}u_1 + vv_1),\\
  r_{3,t} &= \frac{1}{2} q_2\cos(2w)u + \frac{1}{2}q_2\sin(2w)v_1 \\
  &+ \frac{1}{2}p_2\cos(2w)(uv + v_2) + \frac{1}{2}p_2\sin(2w)( -u_1 + vv_1) \\
  &- q_1\cos(2w)(\frac{1}{2}u_1 + vv_1) + q_1\sin(2w)(uv - \frac{1}{2}v_2) \\
  &+ p_1\cos(2w)(\frac{5}{2}uv_1 + \frac{1}{2}u_1v - v^2v_1 - \frac{1}{2}v_3) \\
  &+ p_1\sin(2w)(uv^2 + \frac{1}{2}u_2 + \frac{1}{2}vv_2 + \frac{5}{2}v_1^2) \\
  &+ q_0\cos(2w)(-\frac{3}{2}u^2 - \frac{1}{2}uv^2 + \frac{1}{2}u_2 + \frac{1}{2}vv_2 + v_1^2) \\
  &+ \frac{1}{2}q_0\sin(2w) (-5uv_1 - u_1v - v^2v_1 + v_3) \\
  &+ p_0\cos(2w)(3u^2v + uv^3 - \frac{5}{2}uv_2 - 3u_1v_1 - \frac{3}{2}u_2v - \frac{3}{2}v^2v_2 \\
  &- 3vv_1^2 + \frac{1}{2}v_4) \\
  &+ p_0\sin(2w)(3uu_1 + 6uvv_1 + \frac{3}{2}u_1v^2 - \frac{1}{2}u_3 +
  v^3v_1 - \frac{3}{2}vv_3 - \frac{5}{2}v_1v_2).
\end{align*}

In this extended setting, equation~\eqref{eq:ell*ell} acquires a new
solution of the form
\begin{equation*}
  \begin{pmatrix}
    F\\
    G
  \end{pmatrix}
  =
  \begin{pmatrix}
    \sum_i(F_i^up_i+F_i^vq_i+F_i^wr_i)\\
    \sum_i(G_i^up_i+G_i^vq_i+G_i^wr_i)
  \end{pmatrix},
\end{equation*}
where
\begin{align*}
  F& = (16uu_1 + 12uvv_1 + 6u_1v^2 - 2u_3 - 6vv_3 - 5v_1v_2)p_0\\
  & + (16u^2 + 12uv^2 - 8u_2 - 17vv_2 - 4v_1^2)p_1\\
  & - 3(4u_1 + 5vv_1)p_2 - (8u + 5 v^2)p_3 + p_5\\
  & + (2uv_1 + 5u_1v - v_3)q_0 + (11uv + 2v^3 - 4v_2)q_1- 5v_1q_2 - 4vq_3\\
  & + 2r_3 ( - 2 \cos(2w) u v + \cos(2w) v_2 - \sin(2w) u_1 - 2 \sin(2w) v v_1) \\
  & + 2r_2 (\cos(2w) u_1 + 2 \cos(2w) v v_1 - 2 \sin(2w) u v +\sin(2w) v_2) \\
  & - 3r_1 u_1, \\[1ex]
  G& = (9uv_1 + 6u_1v + 6v^2v_1 - 2v_3)p_0 + (11uv + 2v^3 - 6v_2)p_1 \\
  & - 7v_1p_2 - 4vp_3 + (u_1 + 5vv_1)q_0 + (2u + 5v^2)q_1 - q_3\\
  & - 3v_1r_1 + 2(-\cos(2w) v_1 + \sin(2w) u)r_2 \\
  & + 2(\cos(2w) u + \sin(2w) v_1)r_3.
\end{align*}

In the conventional matrix-operator form this solution looks as
follows
\begin{equation*}
  A'=L+N,
\end{equation*}
where
\begin{equation*}
  L=\begin{pmatrix}
    L_{11}&L_{12}\\
    L_{21}&L_{22}
  \end{pmatrix},\quad
  N=\begin{pmatrix}
    N_{11}&N_{12}\\
    N_{21}&N_{22}
  \end{pmatrix}
\end{equation*}
are $2\times 2$-matrix operators corresponding to the local and
nonlocal parts of $A'$, respectively, and having the following
entries:
\begin{align*}
  L_{11}&=D_x^5 - (8u + 5 v^2)D_x^3 - 3(4u_1 + 5vv_1)D_x^2 \\
  &+ (16u^2 + 12uv^2 - 8u_2 - 17vv_2 - 4v_1^2)D_x \\
  &+16uu_1 + 12uvv_1 + 6u_1v^2 - 2u_3 - 6vv_3 - 5v_1v_2,\\
  L_{12}&=- 4vD_x^3 - 5v_1D_x^2 + (11uv + 2v^3 - 4v_2)D_x +2uv_1 + 5u_1v - v_3,\\
  L_{21}&= - 4vD_x^3 - 7v_1D_x^2 + (11uv + 2v^3 - 6v_2)D_x \\
  &+ 9 u v_1 + 6 u_1 v + 6 v^2 v_1 - 2 v_3,\\
  L_{22}&= - D_x^3 + (2 u + 5 v^2) D_x +(u_1 + 5 v v_1)
  \intertext{and}
  N_{11}&=-3 Y_{1,0}^uD_x^{-1}\circ Y_{1,0}^u-4 Y_{1,1}^u D_x^{-1}\circ Y_{1,1}^u-4 Y_{1,2}^u D_x^{-1}\circ Y_{1,2}^u,\\
  N_{12}&=-3 Y_{1,0}^uD_x^{-1}\circ Y_{1,0}^v-4 Y_{1,1}^u D_x^{-1}\circ Y_{1,1}^v-4 Y_{1,2}^u D_x^{-1}\circ Y_{1,2}^v,\\
  N_{21}&=-3 Y_{1,0}^vD_x^{-1}\circ Y_{1,0}^u-4 Y_{1,1}^v D_x^{-1}\circ Y_{1,1}^u-4 Y_{1,2}^v D_x^{-1}\circ Y_{1,2}^u,\\
  N_{22}&=-3 Y_{1,0}^vD_x^{-1}\circ Y_{1,0}^v-4 Y_{1,1}^v
  D_x^{-1}\circ Y_{1,1}^v-4 Y_{1,2}^v D_x^{-1}\circ Y_{1,2}^v,
\end{align*}
whereas
\begin{equation*}
  Y_{1,0}=\begin{pmatrix}
    Y_{1,0}^u\\
    Y_{1,0}^v
  \end{pmatrix},\quad
  Y_{1,1}=\begin{pmatrix}
    Y_{1,1}^u\\
    Y_{1,1}^v
  \end{pmatrix},\quad
  Y_{1,2}=\begin{pmatrix}
    Y_{1,2}^u\\
    Y_{1,2}^v
  \end{pmatrix},
\end{equation*}
are \emph{symmetries} of the coupled KdV-mKdV system
(see~\cite{KerstenKrasilshchik:CInCKS}) presented in the form
\begin{align*}
  Y_{1,0}^u&=u_1,\\
  Y_{1,0}^v&=v_1;\\[1ex]
  Y_{1,1}^u&=-\cos(2w)(\frac{1}{2}u_1 + vv_1) + \sin(2w)(uv - \frac{1}{2}v_2),\\
  Y_{1,1}^v&=\frac{1}{2}\cos(2w)v_1 - \frac{1}{2}\sin(2w)u;\\[1ex]
  Y_{1,2}^u&=\cos(2w)(uv - \frac{1}{2}v_2) + \sin(2w)(\frac{1}{2}u_1 + vv_1),\\
  Y_{1,2}^v&= - \frac{1}{2}\cos(2w)u - \frac{1}{2}\sin(2w)v_1 .
\end{align*}

\begin{remark}
  Expressions for the entries of the operator $N$ were obtained as it
  was indicated in Remark~\ref{rem:correspondence}.
\end{remark}

As above, simulating the techniques developed for the local theory, we
have checked that
\begin{equation*}
  [\![A',A']\!]=0\quad\text{ and }\quad [\![A',A]\!]=0,
\end{equation*}
i.e., $A$ and $A'$ are compatible.

\begin{remark}\label{rem:hierarchies}
  Though system~\eqref{eq:kdv-mkdv} is Hamiltonian with respect to
  $A'$, there does not exist the corresponding Hamiltonian. This can
  be proved using the same techniques we used in the proof of
  Theorem~\ref{thm:Ham-kdv-mkdv}.  Nevertheless, the following facts
  are valid. Recall that~\eqref{eq:kdv-mkdv} possesses a
  \emph{recursion operator}~\cite{KerstenKrasilshchik:CInCKS}. Denote
  this operator by $R$ and note that our Hamiltonian operators are
  related to each other by means of this recursion operator, i.e.,
  \begin{equation*}
    A'=R\circ A.
  \end{equation*}
  Moreover, in the same way one can construct a whole hierarchy of
  pairwise compatible Hamiltonian structures. On the other hand, $R$
  generates a hierarchy of equations in which~\eqref{eq:kdv-mkdv} is the
  first term. Then $A'$ is a Hamiltonian structure for all other
  equations of this hierarchy and these equations possess Hamiltonians
  with respect to $A'$.
\end{remark}

To make our exposition self-contained, we describe the recursion
operator $R$ in the Appendix.

\section*{Appendix: Recursion operator for the coupled KdV-mKdV
system~\cite{KerstenKrasilshchik:CInCKS}}

To construct the recursion operator, it needs to extend the nonlocal
setting introduced above. Namely, we add three new nonlocal variables
$w_1$, $w_{11}$, and $w_{12}$ defined by
\begin{align*}
  w_{1,x}&=u,\\
  w_{1,t}&=3u^2+3uv-u_2-3vv_2;\\[1ex]
  w_{11,x}&=\cos(2w)w_1v+\sin(2w)v^2,\\
  w_{11,t}&=\cos(2w)(3wuv+wv^3-wv_2+uv_1-u_1v-v^2v_1)\\
  &+\sin(2w)(4uv^2+v^4-vv_2-v_1^2);\\[1ex]
  w_{12,x}&=\cos(2w)v^2-\sin(2w)wv,\\
  w_{12,t}&=\cos(2w)(4uv^2-v^4-2vv_2+v_1^2)\\
  &+\sin(2w)(-3wuv-wv^3+wv_2-uv_1+v^2v_1)
\end{align*}
(see Example~\ref{ex:cover}).  Then $R$ is a $2\times 2$-matrix
operator,
\begin{equation*}
  R=\begin{pmatrix}
    R_{11}&R_{12}\\
    R_{21}&R_{22}
  \end{pmatrix}
\end{equation*}
with the entries
\begin{align*}
  R_{11} &=D_x^2 - 4 u - v^2 \\
  &- Y_{1,1}^u D_x^{-1}\circ \psi_{1,2}^u  + Y_{1,2}^u D_x^{-1}\circ \psi_{1,1}^u - \frac{3}{2} Y_{1,0}^u D_x^{-1}, \\
  R_{12}  &= 2 v D_x^2 + v_1 D_x - 3 u v + 2 v_2  \\
  & - Y_{1,1}^u D_x^{-1}\circ \psi_{1,2}^v + Y_{1,2}^u D_x^{-1}\circ \psi_{1,1}^v + Y_{2,1}^u D_x^{-1},  \\
  R_{21} &= - 2 v  \\
  & - Y_{1,1}^v D_x^{-1}\circ \psi_{1,2}^u   + Y_{1,2}^v D_x^{-1}\circ \psi_{1,1}^u - \frac{3}{2} Y_{1,0}^v D_x^{-1},  \\
  R_{22}  &= D_x^2  - 2 u - v^2 \\
  & - Y_{1,1}^v D_x^{-1}\circ \psi_{1,2}^v + Y_{1,2}^v D_x^{-1}\circ \psi_{1,1}^v + Y_{2,1}^v D_x^{-1},  \\
\end{align*}
where $Y_{1,0}$, $Y_{1,1}$, $Y_{1,2}$ are the same symmetries that
enter the expression for the nonlocal Hamiltonian structure $A'$,
$Y_{2,1}$ is another symmetry with the components
\begin{align*}
  Y_{2,1}^u &=\cos(2 w) ( - w_{11} u_1 - 2 w_{11} v v_1 + 2 w_{12} u v - w_{12} v_2) \\
  &+ \sin(2 w) (2 w_{11} u v - w_{11} v_2 + w_{12} u_1 + 2 w_{12} v v_1) \\
  &- 2 u v_1 - 3 u_1 v - 2 v^2 v_1 + v_3,\\
  Y_{2,1}^v&=\cos(2 w) (w_{11} v_1 - w_{12} u) - \sin(2 w) (w_{11} u + w_{12} v_1) \\
  &- (u_1 + v v_1),
\end{align*}
while
\begin{align*}
  \psi_{1,1}^u&=- \sin(2 w), \\
  \psi_{1,1}^v&=2 (2 \sin(2 w) v - w_{12}), \\[1ex]
  \psi_{1,2}^u&=- \cos(2 w), \\
  \psi_{1,2}^v&=2 (2\cos(2 w) v + w_{11}).
\end{align*}
Note that $\psi_{1,1}=(\psi_{1,1}^u,\psi_{1,1}^v)$ and
$\psi_{1,2}=(\psi_{1,2}^u,\psi_{1,2}^v)$ are the generating functions
for (nonlocal) conservation laws.

\section*{Acknowledgments}
I.~K. and A.~V. are grateful to the University of Twente, where this
research was done, for hospitality.  The work of A.~V. was supported
in part by the NWO and FOM (The Netherlands).  We are also grateful to
our first reader, S.~Igonin, and to the referee for their useful
remarks.

\end{document}